\documentclass[%envcountsame,
    envcountsect,
    envcountsame,
                ]{svmult}

\usepackage{type1cm}        % activate if the above 3 fonts are
\usepackage{tikz}
\usetikzlibrary{matrix}
\usepackage{xcolor}
\usepackage{comment}

\usepackage{cite}
\usepackage{graphicx}        % standard LaTeX graphics tool
\usepackage{multicol}        % used for the two-column index
\usepackage[bottom]{footmisc}% places footnotes at page bottom

\usepackage{enumitem} %enumeration
\usepackage{mathtools}
\usepackage{algorithm}
\usepackage[noend]{algpseudocode}
\usepackage{mleftright}
\mleftright

\makeatletter
\def\BState{\State\hskip-\ALG@thistlm}
\makeatother

\usepackage{placeins}
\usepackage{newtxtext}       %
\usepackage{newtxmath}       % selects Times Roman as basic font

\usepackage{jack}
\usepackage{tikzimages}

 \usepackage[nodayofweek]{datetime}
 
\usepackage{pifont}% http://ctan.org/pkg/pifont

\usepackage{
            ulem		% use \sout{...}
} 

\usepackage{hyperref}

\newcommand{\ootimes}{\,\underline{\otimes}\,}

\newcommand{\dlpar}{\left(\mkern-3mu\left(}
\newcommand{\drpar}{\right)\mkern-3mu\right)}

\makeindex             % used for the subject index
\usepackage{listings}

\begin{document}

\title*{An introduction to tensors for path signatures}
\titlerunning{An introduction to tensors for path signatures}
% Use \titlerunning{Short Title} for an abbreviated version of
% your contribution title if the original one is too long
\author{Jack Beda, Gon\c calo dos Reis %\orcidID{0000-0002-4993-2672}
and Nikolas Tapia}
% Use \authorrunning{Short Title} for an abbreviated version of
% your contribution title if the original one is too long
\institute{Jack Beda \at The University of Edinburgh, James Clerk Maxwell Building, Peter Guthrie Tait Rd, Edinburgh EH9 3FD, UK \email{jack@beda.ca}
\and
Gon\c calo dos Reis \at Maxwell Institute for Mathematical Sciences, School of Mathematics, University of Edinburgh, JCMB, Peter Guthrie Tait Rd, Edinburgh EH9 3FD, UK \email{G.dosReis@ed.ac.uk}
\and
Nikolas Tapia \at Weierstrass Institute, Mohrenstr.~39, 10117 Berlin, Germany, \email{tapia@wias-berlin.de}
}
\maketitle

%%%%%%%%%%%%%%%%%%%%%%%%%%%%%%%%%%%%%%%%%%%%%%%%
%%%%%%%%%%%%%%%%%%%%%%%%%%%%%%%%%%%%%%%%%%%%%%%%%%%%
%%%%%%%  BEGIN ABSTRACTS
%%%%%%%%%%%%%%%%%%%%%%%%%%%%%%%%%%%%%%%%%%%%%%%%%%%%
% %% INSTRUCTIONS
% Please use the ’starred’ version of the new Springer abstract command for
% typesetting the text of the online abstracts (cf. source file of this chapter template abstract) and include them with the source files of your manuscript. Use the plain abstract command if the abstract is also to appear in the printed version of the book

%% THIS APPEARS IN THE BOOK CHAPTER
\abstract{We present a fit-for-purpose introduction to tensors and their operations. 
It is envisaged to help the reader become acquainted with its underpinning concepts for the study of path signatures. 
The text includes exercises, solutions and many intuitive explanations. 
The material discusses direct sums and tensor products as two possible operations that make the Cartesian product of vectors spaces a vector space. The difference lies in linear Vs.~multilinear structures -- the latter being the suitable one to deal with path signatures. The presentation is offered to understand tensors in a deeper sense than just a multidimensional array. 
The text concludes with the prime example of an algebra in relation to path signatures: the \textit{tensor algebra}.\newline \newline
This manuscript is the extended version (with two extra sections) of a chapter to appear in Open Access in a forthcoming Springer volume ``\textit{Signatures Methods in Finance: An Introduction with Computational Applications}". The two additional sections here discuss the factoring of tensor product expressions to a minimal number of terms. This problem is relevant for the path signatures theory but not necessary for what is presented in the book. Tensor factorization is an elegant way of becoming familiar with the language of  tensors and tensor products. 
A GitHub repository is attached.
}

%% THIS APPEARS for the online part 
% \abstract*{
% }

%%%%%%%%%%%%%%%%%%%%%%%%%%%%%%%%%%%%%%%%%%%%%%%%%
%%%%%%%  END ABSTRACTS
%%%%%%%%%%%%%%%%%%%%%%%%%%%%%%%%%%%%%%%%%%%%%%%%%%%%
%%%%%%%%%%%%%%%%%%%%%%%%%%%%%%%%%%%%%%%%%%%%%%%%%%%%%%%%%%%%%%

\section{Introduction}

Throughout mathematics, computer science, and physics, the term \emph{tensor} is used to describe a myriad of similar, but fundamentally different mathematical objects. For amusement, we invite the reader to visit the ``100 Questions: A Mathematical Conventions Survey'' and check\footnote{Find the survey here:  
\url{https://cims.nyu.edu/~tjl8195/survey/results.html\#q45}, at the time this manuscript was written.} ``Question 45: What is a tensor?" That answer's wide distribution seems to hint at a gap, as folks knowledge goes, to the mathematical meaning of a tensor (even among the informed?). 
It is thus necessary to be clear and unambiguous with our definitions and language so that at the end of this chapter the reader will be able to agree with us on the answer. For practical purposes, it often suffices to describe a tensor as a multidimensional array that extends the concept of a matrix \cite{rabanser2017introduction}. 
This is surely true, but only after a certain structure on the underlying space is assumed. Tensors are much more, especially when the underlying spaces are infinite-dimensional \cite{reedbarry1980methodmodernphysics,Diestel2008book,Hackbusch2019book}.

For path signatures, seeing tensors as multidimensional arrays is a good starting point, but their power is only fully realized with an understanding of the operations that go with them: that is, with the algebra mixed in. The assumption behind the next sections is that the reader is largely unfamiliar with tensors (but having heard of multidimensional arrays). 
\smallskip 

\textit{Why do tensors come up in path signatures and how to read this chapter?} We saw in \cite{chevyrev2016primer} that path signatures rely heavily on iterated integrals and involve products of path components across different times. The path signature is a sequence of terms encoding information about a path at different levels of complexity: the first-level signature captures linear information via integrals of individual components; the second-level signature captures pairwise interactions, so-called bilinear relationships, like $\iint\mathrm{d}x^i\mathrm{d}x^j$, meaning they combine two inputs in a way that is linear in each (with $x^i$ fixed, we have linearity in $x^j$, and vice-versa) whereas joint linearity does not hold; higher levels of the signature require more iterated integrals and in turn one requires multilinear relationships to describe these complicated interactions.   
Tensors and multilinear maps are natural and powerful mathematical tools to represent and analyze these interactions.   

Section \ref{firstIntro2Tensors} introduces the basic definitions and properties of tensors by guiding the reader to first understand the difference between linear and bilinear operators, and how linear operators can be recovered from bilinear ones via the tensors product operation and the so-called \textit{universal property}. 
 
Section \ref{section:AlittleBITmore} then goes into a few other deeper properties and connects more explicitly tensors and algebras by introducing the \textit{tensor algebra}. We will leave to subsequent chapters the development of further concepts like shuffle algebras, and power series to obtain exponentials and logarithms (alluded to in \cite{chevyrev2016primer}). % 

\begin{acknowledgement}
NT acknowledges support from DFG CRC/TRR 388 “Rough Analysis, Stochastic Dynamics and Related Fields”, Project B01. 

GdR acknowledges partial support by the Engineering and Physical Sciences Research Council (EPSRC) [grant number EP/R511687/1], from the FCT – Fundação para a Ciência e a Tecnologia, I.P., under the scope of the projects \href{https://doi.org/10.54499/UIDB/00297/2020}{UIDB/00297/2020} and \href{https://doi.org/10.54499/UIDP/00297/2020}{UIDP/00297/2020} (Center for Mathematics and Applications, NOVA Math), and by the UK Research and Innovation (UKRI) under the UK government’s Horizon Europe funding Guarantee [Project APP55638].
\medskip

\textbf{Note}\ 
A GibHub repository related to the content of Sections \ref{sec:TensorRankSection} and \ref{sec:factoringtensorProducts} can be found in 
\begin{center}
\small
\url{https://github.com/jfbeda/tensor_factoring} .
\normalsize
\end{center}

\end{acknowledgement}

%% FOR THE ARXIV VERSION
% \begin{trailer}{Code Notebook}
%     The complete and documented code for this chapter can be found in the following GitHub repository: 
% \begin{center}
% \small
% \url{https://github.com/jfbeda/tensor_factoring} .
% \normalsize
% \end{center}
% \end{trailer}

%%%%%%%%%%%%%%%%%%%%%%%%%%%%%%%%%%%%%%%%%%%%%%%%%%%%%%%%%%%%%%%%
%%%%%%%%%%%%%%%%%%%%%%%%%%%%%%%%%%%%%%%%%%%%%%%%%%%%%%%%%%%%%%%%
%%%%%%%%%%%%%%%%%%%%%%%%%%%%%%%%%%%%%%%%%%%%%%%%%%%%%%%%%%%
%%%%%%%%%%%%%%%%%%%%%%%%%%%%%%%%%%%%%%%%%%%%%%%%%%%%%%%%%%%%%%%%%%%%%
%%%%%%%%%%%%%%%%%%%%%%%%%%%%%%%%%%%%%%%%%%%%%%%%%%%%%%%%%%%%%%%%

\section{A brief introduction to tensors}
\label{firstIntro2Tensors}

This section introduces two operations: \textit{direct sums} and \textit{tensor products}, two different ways of making new vector spaces out of old ones. Formally, each is a way of equipping the Cartesian product of vector spaces, $U\times V$, with a linear structure. The first leads to linear operators while the second leads to bilinear ones. While both are related, in many aspects they behave very differently\footnote{For instance, there is no open mapping theorem for bilinear surjective maps, nor is there a Hahn-Banach theorem for bilinear continuous forms.}. In a nutshell, bilinear maps exhibit separate linearity in $U$ and $V$ while linear maps exhibit global linearity in $U\times V$. 
The distinction is particularly important in the algebra of tensors, where bilinear maps give rise to the tensor product structure through the so-called \textit{universal property}. 
Linear maps correspond to mappings on the direct product space.

\subsection{Direct sums}
Let us recall that for two sets \(X\) and \(Y\), their Cartesian product \(X\times Y\) is defined as
\[
  X\times Y\coloneqq\left\{ (x,y):x\in X,y\in Y \right\},
\]
that is, the set of all ordered pairs where the first is an element of \(X\) and the second an element of \(Y\).

Given vector spaces \(U\) and \(V\), their Cartesian product does not immediately have a linear structure (i.e. is not immediately a vector space).
In other words, after constructing the \emph{set} \(U\times V\) it is not clear how to add two ordered pairs, or multiply them by scalars. We must \emph{define} a way to add and scale the elements of this set, and it turns out there are multiple, sensible and useful definitions. Direct sums are the simplest way to equip the Cartesian product of two (or more) vector spaces with a linear structure of its own.  

\begin{definition}
\label{def:direct.sum}
Let \(U,V\) be vector spaces. On the Cartesian product \(U\times V\) we define the following \textit{linearity} operations: for \(\mathbf{u},\mathbf{u}_1,\mathbf{u}_2\in U\), \(\mathbf{v},\mathbf{v}_1,\mathbf{v}_2\in V\) and \(\lambda\in\R\)
  \begin{align}
    (\mathbf{u}_1,\mathbf{v}_1)+(\mathbf{u}_2,\mathbf{v}_2)&\coloneqq (\mathbf{u}_1+\mathbf{u}_2,\mathbf{v}_1+\mathbf{v}_2),\label{eq:direct.sum.1}
    \\
    \lambda(\mathbf{u},\mathbf{v})&\coloneqq(\lambda\mathbf{u},\lambda\mathbf{v}).\label{eq:direct.sum.2}
  \end{align}
\end{definition}
One can check that \(U\times V\), equipped with the operations from \eqref{eq:direct.sum.1} and \eqref{eq:direct.sum.2}, is a vector space, which it is customary to denote as \(U\oplus V\), i.e the \emph{direct sum} of \(U\) and \(V\). 

\begin{exercise}\label{ex:dsum}
State the axioms that define a \textit{vector space} and show that \(U\oplus V\) is indeed a vector space.
\end{exercise}

It can also be checked that if \(B_U\) and \(B_V\) are bases for \(U\) and \(V\) (respectively) with $0_U$ and $0_V$ denoting the zero elements of $U$ and $V$ (respectively) then the  set
\[
  B = \{(\mathbf{u},0_V) : \mathbf{u}\in B_U\} \cup\{ (0_U,\mathbf{v}):\mathbf{v}\in B_V \}
\]
is a basis for \(U\oplus V\). It follows immediately that \(\dim(U\oplus V)=\dim (U)+\dim (V)\).

\begin{example}
  \label{ex:direct.sum}
  Let \(U=\mathbb{R}^3\) and \(V=\mathcal{M}_{2\times 2}(\mathbb{R})\) be the space of real \(2\)-by-\(2\) matrices.
  A generic element of \(U\oplus V\) is an ordered pair \((\mathbf{u},\mathbf{A})\) for a vector \(\mathbf{u}\in\mathbb{R}^3\) and a matrix \(\mathbf{A}\in\mathcal{M}_{2\times 2}(\mathbb{R})\). A concrete example would be letting 
  \begin{align*}
    (\mathbf{u},\mathbf{A}) = \left(\begin{bmatrix}3\\-1\\4\end{bmatrix}, \begin{bmatrix}1&-1\\2&3\end{bmatrix}\right)
    \qquad \textrm{and}\qquad 
    (\mathbf{v},\mathbf{B}) = \left(\begin{bmatrix}-2\\1\\1\end{bmatrix}, \begin{bmatrix}5&3\\-1&2\end{bmatrix}\right)
  \end{align*}
  then we have
  \[
    (\mathbf{u},\mathbf{A})+(\mathbf{v},\mathbf{B}) = \left( \begin{bmatrix}1\\0\\5\end{bmatrix},\begin{bmatrix}6&2\\1&5\end{bmatrix} \right).
  \]
A basis for \(U\oplus V\) is (\textit{with a slight abuse of notation for what `0' means})
  \[
    \left\{ \left( \begin{bmatrix}1\\0\\0\end{bmatrix},0 \right), \left( \begin{bmatrix}0\\1\\0\end{bmatrix}, 0 \right),\left( \begin{bmatrix}0\\0\\1\end{bmatrix},0 \right), \left( 0, \begin{bmatrix}1&0\\0&0\end{bmatrix} \right), \left( 0, \begin{bmatrix}0&1\\0&0\end{bmatrix} \right), \left( 0, \begin{bmatrix}0&0\\1&0\end{bmatrix} \right), \left( 0, \begin{bmatrix}0&0\\0&1\end{bmatrix} \right)  \right\}.
  \]

  We observe that indeed \(\dim(U\oplus V) = 7= \dim(\R^3)+\dim (\mathcal{M}_{2\times 2}(\mathbb{R}))\).
\end{example}

\begin{remark}[On notation]
  Elements of \(U\oplus V\) can also be written \textit{additively}, that is \(\mathbf{u}+\mathbf{v}\) denotes the vector \((\mathbf{u},\mathbf{v})\in U\oplus V\). This notation is harmless because of Definition \ref{def:direct.sum}, as it behaves in the expected way. We can then restate the two linearity equations \eqref{eq:direct.sum.1} and \eqref{eq:direct.sum.2} in more natural notation:
\[
  \mathbf{u}_1+\mathbf{v}_1+\mathbf{u}_2+\mathbf{v}_2 = \mathbf{u}_1+\mathbf{u}_2+\mathbf{v}_1+\mathbf{v}_2,\qquad \lambda(\mathbf{u}+\mathbf{v}) = \lambda \mathbf{u}+\lambda \mathbf{v},
\]
and we note that there is no confusion in the first equality with the two different meanings of \(+\) since due to the associativity and commutativity of addition, the four terms can be rearranged in an arbitrary way to give the same result\footnote{Formally, there is a canonical isomorphism between \(U\oplus V\) and \(V\oplus U\), so that the pairs \((\mathbf{u},\mathbf{v})\) and \((\mathbf{v},\mathbf{u})\) can be identified.}.  
 Importantly, both the bracket notation, \((\mathbf{u},\mathbf{v})\), and the additive notation, \(\mathbf{u}+\mathbf{v}\), are used interchangeably in the path signature community.
\end{remark}

% \begin{remark}
Definition \ref{def:direct.sum} generalizes easily to a finite number of summands: if \(U_1,\dotsc,U_n\) are vector spaces, the set \(U_1\times\dotsm\times U_n\) carries a linear structure given by componentwise addition and multiplication by scalars. The resulting vector space is denoted by \(U_1\oplus\dotsb\oplus U_n\). 
Extending this concept to infinite families requires some care.
% \end{remark}

\begin{definition}\label{def:inf.dsum}
  Consider an infinite index set \(I\) and let \((U_i:i\in I)\) be a family of vector spaces.
  The direct sum is defined to be the set of all sequences \((\mathbf{u}_i:i\in I)\) such that \(\mathbf{u}_i\neq 0\) for finitely many indices \(i\in I\). Addition and scalar multiplication are defined componentwise. The resulting vector space is denoted by
  \[
    \bigoplus_{i\in I}U_i.
  \]
\end{definition}

  The finiteness constraint in this definition means that direct sum is a subset of the Cartesian product of the spaces, that is,
  \[
    \bigoplus_{i\in I}U_i\subseteq\prod_{i\in I}U_i,
  \]
  where \(\prod_{i\in I}U_i\) is simply the set of all sequences indexed by \(I\).

For our purposes it will be enough to consider countable families of vector spaces, that is, we will take \(I=\mathbb{N}\).
In this case, elements of \(\bigoplus_{n\in\mathbb{N}}U_n\) may be denoted as
\[
  (\mathbf{u}_0,\mathbf{u}_1,\mathbf{u}_2,\dotsc)
\]
with the convention that there is only a finite number of non-zero elements in the sequence. 
The Cartesian product also carries the same linear structure and is indeed a vector space, whose elements consist of arbitrary \(\mathbb{N}\)-indexed sequences, which are still denoted as  
% \[
 $ (\mathbf{u}_0,\mathbf{u}_1,\mathbf{u}_2,\dotsc) $ 
% \]
where all entries may be non-zero.

It should be noted that in the case that \(I\subset\mathbb{N}\) is a finite set, say \(I=\{1,\dotsc,N\}\), both spaces coincide but the inclusion becomes strict as soon as \(I\) is countable. In particular, when dealing with finite collections of spaces there is no ambiguity in the notation.

The main application of direct sums in the world of signatures is to decompose a vector space in terms of subspaces of objects sharing similar ``shape'' properties.

\begin{definition}
\label{def:gvs}
  A vector space \(V\) is said to be \textit{graded} if it can be decomposed as a direct sum:
  \[
    V = \bigoplus_{n\in\N}V_n.
  \]
  The subspace \(V_n\) is called the \emph{homogeneous component of degree \(n\)}. For \(v\in V_n\) we write \(|v|=n\) for its degree.  
\end{definition}
We also note that this definition includes the case of finitely many summands, in which case there is \(N\in\mathbb{N}\) such that \(V_n=\{0\}\) for all \(n>N\).

\subsection{Tensor product and Tensors}
We have now seen how the \textit{direct product} is one way of equipping the Cartesian product of two vector spaces with a vector space structure. There is another way, the \textit{tensor product}, in many ways similar, but with a structure such that it is compatible with multilinear functions in a way to be made precise later.
\begin{definition}\label{def:tensor.prod}
  Let \(U,V\) be vector spaces. On the set \(U\times V\) define the following \textit{bilinearity} operations: for \(\mathbf{u}, \mathbf{u}_1, \mathbf{u}_2\in U\), \(\mathbf{v}, \mathbf{v}_1, \mathbf{v}_2\in V\) and \(\lambda\in\mathbb{R}\)
  \begin{align*}
    (\mathbf{u}_1, \mathbf{v}) + (\mathbf{u}_2, \mathbf{v}) &\coloneqq (\mathbf{u}_1 + \mathbf{u}_2, \mathbf{v}),
    \\
    (\mathbf{u}, \mathbf{v}_1) + (\mathbf{u}, \mathbf{v}_2) &\coloneqq (\mathbf{u}, \mathbf{v}_1 + \mathbf{v}_2),
    \\
   \lambda(\mathbf{u},\mathbf{v})
    &\coloneqq(\lambda\mathbf{u},\mathbf{v}) 
        \quad\textrm{and} 
        \footnotemark
            \quad
   \lambda(\mathbf{u},\mathbf{v})\coloneqq (\mathbf{u},\lambda \mathbf{v}).
      \end{align*}
\end{definition}
\footnotetext{
This last element of the definition is actually \emph{imposing} the vectors \((\lambda\mathbf{u},\mathbf{v})\) and \((\mathbf{u},\lambda\mathbf{v})\) to be equal in \(U\otimes V\). This could be formalized by the use of quotient spaces, but doing such would involve a level of additional complexity not really necessary at the moment. (This also takes care of the apparent non-uniqueness of \(0_{U\otimes V}\) hinted at in Exercise 2.10.)}
As in the direct sum case (Definition \ref{def:direct.sum}), it can be verified that \textit{$U\times V$ equipped with the bilinearity operations} is a vector space. We denote 
the \textit{tensor product} of \(U\) and \(V\) as \(U\otimes V\). For elements of the tensor product space, we write \(\mathbf{u} \otimes \mathbf{v}\coloneqq (\mathbf{u},\mathbf{v})\in U\otimes V\).
By a slight abuse of language we also refer to \(\mathbf{u}\otimes\mathbf{v}\) as the tensor product of the vectors \(\mathbf{u}\in U\) and \(\mathbf{v}\in V\).
As we will see later in Section \ref{section:AlittleBITmore}, this name is justified.

\begin{exercise}
  \label{ex:tensor}
  Show that \(U\otimes V\) is a vector space (recall Exercise \ref{ex:dsum}).
\end{exercise}

\begin{exercise}\label{ex:otimes.familiarity} In this exercise we see  why the \(\otimes\) notation is an intuitive way of writing the tensor product: 
% \begin{enumerate}[label=(\alph*)]
% \item 
(a) Rewrite the bilinearity operations of Definition \ref{def:tensor.prod} using the notation \(\mathbf{u} \otimes \mathbf{v}\) instead of \((\mathbf{u},\mathbf{v})\); 
% \item 
(b) Expand \( (\mathbf{u}_1+\mathbf{u}_2)\otimes(\mathbf{v}_1+\mathbf{v}_2)\), where of course the addition \(\mathbf{u}_1+\mathbf{u}_2\in U\) is simply the addition in \(U\), and the same for \(V\); 
% \item  
(c) We have \(\lambda \mathbf{u}\otimes \lambda \mathbf{v} \propto (\mathbf{u}\otimes \mathbf{v})\). What is the constant of proportionality?
% \end{enumerate}

\end{exercise}
Contrary to the direct sum case, for tensor products its not always possible to write  \(\mathbf{u}_1\otimes \mathbf{v}_1+\mathbf{u}_2\otimes \mathbf{v}_2\) as a single tensor product of a vector in \(U\) with a vector in \(V\)
\footnote{As an aside, it is this property of tensor product spaces that mathematically captures the phenomenon of entanglement of quantum particles. A quantum state like \(|\varphi \rangle = \begin{bmatrix}
    1\\
    0
\end{bmatrix} \otimes \begin{bmatrix}
    0\\
    1
\end{bmatrix} + \begin{bmatrix}
    0\\
    1
\end{bmatrix}\otimes \begin{bmatrix}
    1\\
    0
\end{bmatrix}\) is entangled as it cannot be written as \(\begin{bmatrix}
    a\\
    b
\end{bmatrix}\otimes \begin{bmatrix}
    c\\
    d
\end{bmatrix}\) for any scalars \(a,b,c,d\). This entangled state says ``I have two particles, and their spins are always opposite, but I cannot know which one is spin-up, and which one is spin-down".}.

%\(|\varphi \rangle = |\!\! \uparrow \rangle \otimes |\!\!\downarrow\rangle + |\!\! \downarrow \rangle \otimes |\!\!\uparrow\rangle\)

\begin{exercise}\label{ex:tensor.0}
  Show that \(0_U\otimes \mathbf{v}=\mathbf{u}\otimes 0_V=0_{U\otimes V}\) for all \(\mathbf{u}\in U, \mathbf{v}\in V\).
\end{exercise}

In the same vein as the comment offered just before Definition \ref{def:inf.dsum}, the Definition \ref{def:tensor.prod} admits a straightforward generalization to finitely many vector spaces and extending this concept to infinite families requires some care.  
In particular, we write
\[
    U^{\otimes n} \coloneqq \underbrace{U\otimes U\otimes\dotsm\otimes U}_{n\text{ times }}.
\]

\begin{proposition}\label{prp:tensor.basis}
  For \(U,V\) vector spaces with bases \(B_U,B_V\), respectively, the set
  \[
    B = \{\mathbf{u}\otimes \mathbf{v}:\mathbf{u}\in B_U,\mathbf{v}\in B_V\}
  \]
  is a basis for \(U\otimes V\).
  It follows that \(\dim(U\otimes V)=\dim (U) \cdot \dim (V)\).
\end{proposition}

We can tell immediately by the dimension of the spaces that the direct product and tensor product produce fundamentally different vector structures on the same set\footnote{Recall \(\dim(U\oplus V) = \dim(U) + \dim(V)\), whereas \(\dim(U\otimes V) = \dim(U) \cdot \dim(V)\).}. What distinguishes the tensor product vector space from all other possible linear structures is the following \emph{universal property}.
\begin{theorem}
\label{theo:UniversalityProperty}
  Let \(U,V,W\) be vector spaces, and let \(f\colon U\times V\to W\) be a bilinear map. That is, $f$ satisfies for all \(\mathbf{u}, \mathbf{u}_1, \mathbf{u}_2\in U\), \(\mathbf{v}, \mathbf{v}_1, \mathbf{v}_2\in V\) and \(\lambda\in\R\)
  \begin{gather*}
    f(\mathbf{u}_1+\mathbf{u}_2,\mathbf{v}) = f(\mathbf{u}_1, \mathbf{v}) + f(\mathbf{u}_2, \mathbf{v}),
    \qquad f(\mathbf{u}, \mathbf{v}_1+\mathbf{v}_2) = f(\mathbf{u}, \mathbf{v}_1)+f(\mathbf{u}, \mathbf{v}_2),
    \\
    \textrm{and}\quad  f(\lambda \mathbf{u}, \mathbf{v}) = f(\mathbf{u}, \lambda \mathbf{v}) = \lambda f(\mathbf{u},\mathbf{v}).
  \end{gather*}
  Then, there exists a unique \textbf{linear} function \(\hat f\colon U\otimes V\to W\) such that \(f(\mathbf{u},\mathbf{v}) = \hat f(\mathbf{u}\otimes \mathbf{v})\).
\end{theorem}
We say that bilinear functions \emph{factor through} the tensor product.
In fact, the tensor product is characterized by this property, in the sense that any other vector space \(Z\) equipped with a map \(\ootimes\colon U\times V\to Z\) factorizing bilinear functions must be isomorphic to \(U\otimes V\). 
In other words, the tensor product is the unique (up to isomorphism) vector space having this property.
For the sake of simplicity we will omit the proof of this result, but refer the interested reader to the classical texts \cite{Diestel2008book,Hackbusch2019book}. 
At first sight bilinear maps are not quite as powerful as linear maps, nonetheless, the universal property fixes things as it allows to write the bilinear map as a linear map of the tensor product and thus recovers the neat results of linear maps that were not available (at the cost of using tensor products). 

The next example highlights the difference between direct sums and tensor products.
\begin{example}[Direct sums $\oplus$ Vs tensor products $\otimes$]
\label{ex:tensor.sum}
  Let \(U=V=\mathbb{R}\). 
  The direct sum \(U\oplus V\) satisfies \(U\oplus V\cong\mathbb{R}^2\). Indeed, elements of \(U\oplus V\) are ordered pairs of real numbers with component-wise addition and scalar multiplication. Moreover, a basis for \(U\oplus V\) is \(\{(1,0),(0,1)\}\) which is the canonical basis of \(\mathbb{R}^2\), so \(\dim(U\oplus V)=2\).

  On the contrary, we will show now that \(U\otimes V\) satisfies  \(U\otimes V\cong \mathbb{R}\) which is \textit{obviously not} $\mathbb{R}^2 \cong U\oplus V$. 
  Consider the map \(\varphi\colon U\times V\to\mathbb{R}\) given by \(\varphi(x,y)=xy\); this map $\varphi(x,y)$ is clearly bilinear.
  By the universal property, there exists a unique map \(\hat{\varphi}\colon U\otimes V\to\mathbb{R}\), given by \(\hat{\varphi}(x\otimes y)=xy\).
  The map \(\varphi\) is injective since the equation \(\hat{\varphi}(x\otimes y)=0\) implies that either \(x=0\) or \(y=0\); in any case \(x\otimes y=0\) by Exercise \ref{ex:tensor.0}.
  Finally, if \(\lambda\in\mathbb{R}\) then \(\hat{\varphi}(\lambda\otimes 1)=\lambda\) so that \(\hat{\varphi}\) is surjective.
  In particular, \(\mathbb{R}\otimes\mathbb{R}\) is spanned by the vector \(1\otimes 1\) so that \(\dim(U\otimes V)=1\).
\end{example}

\begin{exercise}\label{ex:tensor.R}
  Show that if \(U\) is any vector space, then \(U\otimes\mathbb{R}\) and \(\mathbb{R}\otimes U\) are isomorphic to \(U\). (\emph{Hint: generalize Example \ref{ex:tensor.sum}.})
\end{exercise}

The word \textit{tensor} has many different meanings across different fields \cite{reedbarry1980methodmodernphysics,Diestel2008book,Hackbusch2019book}.  
We will mostly be interested in the case where we are taking tensor products of a finite number of finite-dimensional vector spaces, represented as \(\mathbb{R}^d\) for some integer \(d\ge 1\).
In this setting, tensors may be represented in a simpler, more concrete way, by working with canonical bases. 
It is at at this specific juncture (assuming a basis) that it is intuitive to define a tensor as a multidimensional array -- 
see Remark \ref{rem:tensor.coord}.

\begin{definition}[Tensor: order and shape] 
\label{defn: tensor}  
Take \(n\ge 1\), set \(d_1,\dotsc,d_n\ge 1\) and take the associated vector spaces $\R^{d_j}$ for $j=1,\dotsc,n$.  
An \emph{order \(n\) tensor} of \emph{shape} \((d_1,\dotsc,d_n)\) is an element of the tensor product \(\R^{d_1}\otimes\dotsm\otimes\R^{d_n}\). 
\end{definition}
It should be clear that elements of, say, \(\R^{2}\otimes\R^{4}\), \(\R^{4}\otimes\R^{2}\), and \(\R^3 \otimes \R^3\), are all order 2 tensors, but their \textit{shapes} are all very different -- and, just like matrices, they cannot be added together as the shapes do not match.

\begin{remark}[Basis, vectors and their components, and some notation]
\label{rem:tensor.coord}
  Denote by \(\{\mathbf{e}_1,\dotsc,\mathbf{e}_{d}\}\) the canonical basis of \(\R^{d}\) for some $d\geq 1$. We have seen that the set
\[
  \left\{\mathbf{e}_{i_1}\otimes\dotsm\otimes \mathbf{e}_{i_n}:i_j\in \{1,\dotsc,d_j\}\text{ for all \(j=1,\dotsc,n\)}\right\}
\]
is a basis for \(\R^{d_1}\otimes\dotsm\otimes\R^{d_n}\), which we call the canonical basis.
In particular, an order \(n\) tensor \(\mathbf{T}\) is determined by the \(n\)-dimensional array of its coefficients in this basis:
\[
  \mathbf{T} = \sum_{i_1=1}^{d_1} \dotsb \sum_{i_n=1}^{d_n} \mathbf{T}^{i_1 \dotsb i_n}\mathbf{e}_{i_1}\otimes\dotsm\otimes \mathbf{e}_{i_n}.
\]

As long as we keep this in mind, the assignment \(\mathbf{T}\mapsto (\mathbf{T}^{i_1 \dotsb i_n})\) defines a one-to-one correspondence between elements of the tensor product \(\mathbf{T}\in\R^{d_1}\otimes\dotsm\otimes\R^{d_n}\) and multidimensional arrays \((\mathbf{T}^{i_1\dotsm i_n})\in\R^{d_1\times\dotsm\times d_n}\).
\textit{We remark once again that this isomorphism depends on the fixing of a basis and is, in general, not canonical.}

Thus, \textbf{notation wise}, we shall refer to tensors using either a symbol (i.e. \(\mathbf{T}\)) or in component notation (i.e. \(\mathbf{T}^{ijk}\) -- using superscript notation for its components); for multiple tensors or vectors we use subscript notation, i.e., \(\mathbf{T}_1,\mathbf{T}_2,\dotsc\) or \(\mathbf{u}_1,\mathbf{u}_2,\dotsc\). 

For example, \(\mathbf{C}^{ij}\), \(\mathbf{T}^{ijk}\) and \(\mathbf{Q}^{ijkp}\) refer to the components of tensors \(\mathbf{C}\), \(\mathbf{T}\), and \(\mathbf{Q}\) of order 2, 3, and 4 respectively. 
In particular, the order 1 tensors \(\mathbf{e}_1,\dotsc,\mathbf{e}_d\) constitute the canonical basis of \(\R^d\). The \(i\)th basis vector \(\mathbf{e}_i\) has components given in the canonical basis by
\[
  \mathbf{e}_i^j=\begin{dcases}1&j=i\\0&\text{else}\end{dcases}
  \qquad \textrm{for \(j=1,\dotsc,d\).}
\]
\end{remark}

\begin{example}
\label{example:tensor-matrices-example}
    We see that tensors of order \(1\) and \(2\) can be identified with column vectors and matrices, respectively. A scalar is by convention an order \(0\) tensor.
    For example, \(\mathbf{u}\), \(\mathbf{A}\), and \(\mathbf{T}\) are tensors of order \(1, \,2,\) and \(3\) respectively:
    \begin{equation}
      \begin{gathered}
        \mathbf{u} \coloneqq \begin{bmatrix}
          1\\
          2\\
          3
        \end{bmatrix}\in \R^3,
        \quad \mathbf{A} \coloneqq\begin{bmatrix}
          1 & 2\\
          3 & 4
        \end{bmatrix}\in \R^{2}\otimes\R^2 = (\R^2)^{\otimes 2} \cong\R^{2\times 2},\\\mathbf{T} \coloneqq
      \resizebox{0.25\linewidth}{!}{\tensorone}
        \in(\R^2)^{\otimes 3}\cong\R^{2\times 2\times 2}.
      \end{gathered}
    \end{equation} 
\end{example}

\begin{remark}
  The tensor is an \textit{intrinsic} object, in the sense that tensors do not depend on any choice of basis.
  From Remark \ref{rem:tensor.coord}, we see that a tensor can be uniquely described by a multidimensional array of numbers, but this is only true once we have fixed a basis. For many applications, it is possible, and practical, to think of tensors only in the terms of their components in a particular (e.g. the canonical) basis. 
  Nonetheless, we encourage the reader to be mindful with the language and recognize when they are being loose with the concepts.   
  This is exactly the same idea to how we can think of linear transformations from \(\R^n\) to \(\R^n\) in terms of an \(n\)-by-\(n\) square matrix, once we have fixed a basis. For example, the linear map \((x_1,x_2)\mapsto(x_1+x_2,x_1-x_2)\) can be represented by the matrix \(\begin{pmatrix}1&1\\1&-1\end{pmatrix}\) in the standard basis, but in the eigenbasis, it is represented by the diagonal matrix \(\begin{pmatrix}\sqrt{2}&0\\0&-\sqrt2\end{pmatrix}\).
  The characterizing property of tensors is that given a change of basis, we immediately know how the coordinate representation transforms.
  In this example we can go from the first to the second matrix representation by multiplication with an invertible matrix (diagonalization).
\end{remark}

Writing explicit examples for the tensor product quickly becomes cumbersome, nevertheless we offer a few simple examples.
  \begin{example}[Tensor multiplication]
  \label{example: tensor multiplication}
    Let \(\mathbf{u}\in\R^2\), \(\mathbf{v}\in\R^3\) be vectors, that is, order 1 tensors of shapes \((2)\) and \((3)\), respectively.
    By definition, their tensor product is an order 2 tensor of shape \((2,3)\): \(\mathbf{u}\otimes\mathbf{v}\in\R^2\otimes\R^3\).
    In the canonical basis they are represented by 2 and 3 coefficients, respectively. Namely
    \[
      \mathbf{u}=\mathbf{u}^1\mathbf{e}_1+\mathbf{u}^2\mathbf{e}_2,\quad\mathbf{v}=\mathbf{v}^1\mathbf{e}_1+\mathbf{v}^2\mathbf{e}_2+\mathbf{v}^3\mathbf{e}_3,
    \]
    or in the more traditional column vector notation,
    \[
      \mathbf{u}=\mathbf{u}^1\begin{bmatrix}1\\0\end{bmatrix}+\mathbf{u}^2\begin{bmatrix}0\\1\end{bmatrix}=\begin{bmatrix}\mathbf{u}^1\\\mathbf{u}^2\end{bmatrix}\quad\text{ and }\quad\mathbf{v}=\begin{bmatrix}
      \mathbf{v}^1\\\mathbf{v}^2\\\mathbf{v}^3\end{bmatrix}.
    \]
    By using the bilinearity of the tensor product (Definition \ref{def:tensor.prod}) we may obtain the coordinates of \(\mathbf{u}\otimes\mathbf{v}\) in the canonical basis.
    Indeed, recalling that the components of both vectors are scalars, 
    \begin{equation}\label{eq:uv}
      \begin{split}
        \mathbf{u}\otimes\mathbf{v} &= \left( \mathbf{u}^1\mathbf{e}_1+\mathbf{u}^2\mathbf{e}_2 \right)\otimes\left( \mathbf{v}^1\mathbf{e}_1+\mathbf{v}^2\mathbf{e}_2+\mathbf{v}^3\mathbf{e}_3 \right)\\
                                    &= \begin{multlined}[t]\mathbf{u}^1\mathbf{v}^1\mathbf{e}_1\otimes\mathbf{e}_1+\mathbf{u}^1\mathbf{v}^2\mathbf{e}_1\otimes\mathbf{e}_2+\mathbf{u}^1\mathbf{v}^3\mathbf{e}_1\otimes\mathbf{e}_3\\+\mathbf{u}^2\mathbf{v}^1\mathbf{e}_2\otimes\mathbf{e}_1+\mathbf{u}^2\mathbf{v}^2\mathbf{e}_2\otimes\mathbf{e}_2+\mathbf{u}^2\mathbf{v}^3\mathbf{e}_2\otimes\mathbf{e}_3.
                                    \end{multlined}
      \end{split}
    \end{equation}

    Thus, in the canonical basis the order 2 tensor \(\mathbf{u}\otimes\mathbf{v}\) has components given by \((\mathbf{u}\otimes\mathbf{v})^{ij}=\mathbf{u}^i\mathbf{v}^j\).
    Identifying the canonical basis of order 2 tensors with matrices (i.e., $\mathbf{e}_i \otimes \mathbf{e}_j$ forming the standard \(2\)-by-\(2\) matrix basis), we may write
    \[
      \mathbf{u}\otimes\mathbf{v} = 
      \begin{bmatrix}
        \mathbf{u}^1\mathbf{v}^1\ 
        & 
        \mathbf{u}^1\mathbf{v}^2\ 
        &
        \mathbf{u}^1\mathbf{v}^3 
        \\
        \mathbf{u}^2\mathbf{v}^1\ 
        & 
        \mathbf{u}^2\mathbf{v}^2\ 
        & \mathbf{u}^2\mathbf{v}^3 
      \end{bmatrix}
      \quad \textrm{and likewise}\quad 
    % \]
    % Likewise
    % \[
      \mathbf{v}\otimes\mathbf{u} = 
      \begin{bmatrix}
      \mathbf{v}^1\mathbf{u}^1\ 
      & \mathbf{v}^1\mathbf{u}^2
      \\
      \mathbf{v}^2\mathbf{u}^1
      &\mathbf{v}^2\mathbf{u}^2 \\\mathbf{v}^3\mathbf{u}^1&\mathbf{v}^3\mathbf{u}^2\end{bmatrix}\in\R^3\otimes\R^2.
    \]

    Consider now the matrix 
    \[
      \mathbf{A}\coloneqq\begin{bmatrix}\mathbf{A}^{11}\ &\mathbf{A}^{12}\\\mathbf{A}^{21}&\mathbf{A}^{22}\end{bmatrix}\in\R^2\otimes\R^2.
    \]
    Where, as before, the entries in the matrix notation corresponds to the coordinates in the canonical order 2 tensor basis: 
    \(\mathbf{A}=
    \mathbf{A}^{11}\mathbf{e}_1\otimes\mathbf{e}_1
    + \mathbf{A}^{12}\mathbf{e}_1\otimes\mathbf{e}_2
    + \mathbf{A}^{21}\mathbf{e}_2\otimes\mathbf{e}_1
    +\mathbf{A}^{22}\mathbf{e}_2\otimes\mathbf{e}_2
     \).
    Note that although \(\mathbf{A}\) is an element of \(\R^2\otimes \R^2\), it does not necessarily mean that it can be written as \(\mathbf{u}\otimes\mathbf{v}\) for some \(\mathbf{u},\mathbf{v}\in\R^2\).
    That is, the components \(\mathbf{A}^{ij}\) are not necessarily of the form \(\mathbf{A}^{ij}=\mathbf{u}^i\mathbf{v}^j\) for some vectors \(\mathbf{u}\) and \(\mathbf{v}\) of the appropriate dimensions. 
    In the cases where it \textit{is} possible to find such a decomposition, we say that \(\mathbf{A}\) is a \textit{rank 1 tensor}. On the other hand, $\mathbf{A}$ can always be written as a linear combination of sums of tensor products of some \(\mathbf{u}_i,\mathbf{v}_i\in\R^2\), that is, sums of rank 1 tensors. 

    Continuing with the example, we may compute
    \[\mathbf{u}\otimes \mathbf{A}=
\begin{tikzpicture}[baseline={([yshift=-1.7cm]current bounding box.north)},line width = 0.7pt] 
  \matrix (back) [matrix of math nodes, nodes={draw, minimum size=1cm, anchor=center, fill={rgb, 255:red, 230; green, 225; blue, 225 }}, column sep=-\pgflinewidth, row sep=-\pgflinewidth, xshift=2.3cm, yshift=0.5cm] {
    \mathbf{u}^2\mathbf{A}^{11}& \mathbf{u}^2\mathbf{A}^{12} \\
    \mathbf{u}^2\mathbf{A}^{21} & \mathbf{u}^2\mathbf{A}^{22} \\
  };
  % Define the coordinates for the front-facing matrix
  \matrix (front) [matrix of math nodes, nodes={draw, minimum size=1cm, anchor=center, fill=white, fill opacity=1}, column sep=-\pgflinewidth, row sep=-\pgflinewidth] {
    \mathbf{u}^1\mathbf{A}^{11} & \mathbf{u}^1\mathbf{A}^{12} \\
    \mathbf{u}^1\mathbf{A}^{21} & \mathbf{u}^1\mathbf{A}^{22} \\
  };

  % Draw the connecting lines
  \draw (front-1-1.north west) -- (back-1-1.north west);
  \draw (front-1-2.north east) -- (back-1-2.north east);
  \draw (front-2-1.south west) -- (back-2-1.south west);
  \draw (front-2-2.south east) -- (back-2-2.south east);

  % Draw the edges of the front matrix
  \draw (front-1-1.north west) -- (front-1-2.north east) -- (front-2-2.south east) -- (front-2-1.south west) -- cycle;

  % Draw the edges of the back matrix
  \draw (back-1-1.north west) -- (back-1-2.north east) -- (back-2-2.south east) -- (back-2-1.south west) -- cycle;
\end{tikzpicture}
      \in\R^2\otimes\R^2\otimes\R^2=(\R^2)^{\otimes 3},
    \]
    where the expression on the right-hand side is a matrix-like notation for organizing the components of the order 3 tensor \(\mathbf{u}\otimes\mathbf{A}\). 
  \end{example}
% }

  We will later see, in Section \ref{section:AlittleBITmore}, that the computation performed in equation \eqref{eq:uv} makes use of a larger structure.
  The tensor product can be thought of as a non-commutative analogue of polynomial multiplication, in the sense that it is an associative and bilinear operation. Note that it is, however, not commutative as the results of \(\mathbf{u}\otimes\mathbf{v}\) and \(\mathbf{v}\otimes\mathbf{u}\) are order 2 tensors of different shapes, at least in this example. In general, the results may differ even thought the shapes match.

\begin{example}[Connecting to Signatures of \cite{chevyrev2016primer}
\label{example:basicSigstotenosor}]
  Suppose \(\mathbf{\mathbf{x}}\colon[0,1]\to\R^d\) is a smooth vector-valued path, which simply means, e.g., that for each \(t\in[0,1]\) we may write \(\mathbf{\mathbf{x}}_t=(\mathbf{x}^1_t,\dotsc,\mathbf{x}^d_t)\in\R^d\) in the canonical basis.
  In particular, for each \(t\in[0,1]\), \(\mathbf{x}_t\) is a tensor of order 1.

  We may use the tensor product to compactly write the collection of iterated integrals of \(\mathbf{x}\). That is,
  \[
    \int_0^t\int_0^s\mathrm{d}\mathbf{x}_u\otimes\mathrm{d}\mathbf{x}_s\in\R^d\otimes\R^d
  \]
  is an order 2 tensor (i.e. a \(d\)-by-\(d\) matrix) with components (for $\mathbf{x}$ ``sufficiently nice'')
  \[
    \left( \int_0^t\int_0^s\mathrm{d}\mathbf{x}_u\otimes\mathrm{d}\mathbf{x}_s \right)^{ij} =  \int_0^t\int_0^s\dot{x}^i_u\dot{x}^j_s\,\mathrm{d}u\mathrm{d}s.
  \]
\end{example}

Before delving into tensors properties we can go back to ``Question 45: What is a tensor?'' of the ``100 Questions: A Mathematical Conventions Survey''. We hope to have convinced the reader that a tensor is nothing other than ``an element of a tensor product of vector spaces''.

%%%%%%%%%%%%%%%%%%%%%%%%%%%%%%%%%%%%%%%%%%%%%%%%%%%%%%%%%%%%%%%%
%%%%%%%%%%%%%%%%%%%%%%%%%%%%%%%%%%%%%%%%%%%%%%%%%%%%%%%%%%%%%%%%
%%%%%%%%%%%%%%%%%%%%%%%%%%%%%%%%%%%%%%%%%%%%%%%%%%%%%%%%%%%
%%%%%%%%%%%%%%%%%%%%%%%%%%%%%%%%%%%%%%%%%%%%%%%%%%%%%%%%%%%%%%%%%%%%%
%%%%%%%%%%%%%%%%%%%%%%%%%%%%%%%%%%%%%%%%%%%%%%%%%%%%%%%%%%%%%%%%

\section{A little bit more on tensors: The tensor algebra}
\label{section:AlittleBITmore}

We now venture into a few additional properties of tensors. 
\begin{definition}\label{def:algebra}
  An associative algebra is a vector space \(A\) equipped with a bilinear map \(m\colon A\times A\to A\), called \emph{product}, satisfying the associativity condition
  \[
    m\big(m(\mathbf{x},\mathbf{y}),\mathbf{z}\big)=m\big(\mathbf{x}, m(\mathbf{y},\mathbf{z})\big).
  \]
  The universal property of the tensor product (see Theorem~\ref{theo:UniversalityProperty}) yields that equivalently, it may be represented by a \emph{linear} map \(m\colon A\otimes A\to A\), which is the one we will use from now on.
  
  We say that an algebra \(A\) is \emph{unital} if it has a distinguished element \(1_A\in A\), called the \emph{unit}, satisfying for all \(\mathbf{x}\in A\)
  \[
    m(1_A\otimes \mathbf{x}) = \mathbf{x} = m(\mathbf{x}\otimes 1_A).
  \]
\end{definition}
It is customary to write the product, denoted $\cdot_A$, of two elements of $A$ using infix notation\footnote{Infix notation is a way of writing mathematical (and logic) expressions where operators are placed between the operands they act upon. This is the most familiar notation to us humans and matches how we (humans) interpret math expressions.
%; on the other hand, computers parse expressions more efficiently with prefix or postfix notation.
}, that is, \(\mathbf{x}\cdot_A\mathbf{y}\coloneqq m(\mathbf{x}\otimes\mathbf{y})\). 
Oftentimes, when no confusion can arise, we write simply \(\mathbf{x}\cdot\mathbf{y}\) or even omit the symbol completely and just write \(\mathbf{xy}\) instead.
As we will mostly work with \textit{unital associative algebras}, from here on we simply write  \textit{algebra}. 
In this notation, the condition for \((A,\cdot_A)\) 
to be an algebra can be written as
\[
  (\mathbf{x}\cdot_A \mathbf{y})\cdot_A \mathbf{z} = \mathbf{x}\cdot_A(\mathbf{y}\cdot_A \mathbf{z})
\]
for all \(\mathbf{x},\mathbf{y},\mathbf{z}\in A\) and the unit satisfies \(1\cdot_A \mathbf{x}=\mathbf{x}\cdot_A 1=\mathbf{x}\) for all \(\mathbf{x}\in A\).
We note that bilinearity of the product translates into the distributivity of \(\cdot_A\) over \(+\), e.g.,
\[
  (\mathbf{x}+\mathbf{y})\cdot_A\mathbf{z} = \mathbf{x}\cdot_A\mathbf{z}+\mathbf{y}\cdot_A\mathbf{z},\quad \lambda(\mathbf{x}\cdot_A\mathbf{y})=(\lambda\mathbf{x})\cdot_A\mathbf{y}=\mathbf{x}\cdot_A(\lambda\mathbf{y}),
\]
and so on.
We also remark that we \emph{do not} require the product to be commutative, that is, we do not enforce that \(\mathbf{x}\cdot_A\mathbf{y}=\mathbf{y}\cdot_A\mathbf{x}\) for every \(\mathbf{x},\mathbf{y}\in A\), although this may hold for some pairs of elements.
In case this identity does hold for every \(\mathbf{x},\mathbf{y}\in A\) we say \(A\) is commutative (or Abelian).

\begin{example}
  Let \(\mathcal{M}_{n\times n}(\mathbb{R})\) be the space of \(n\)-by-\(n\) square matrices with real entries with its usual vector space structure (entry-wise addition and scalar multiplication).
  The matrix product \(\mathbf{A}\cdot\mathbf{B}\coloneqq \mathbf{A}\mathbf{B}\) equips \(A\) with the structure of a (non-commutative) associative algebra with unit \(1_A=\mathbf{I}_n\), the \(n\)-by-\(n\) identity matrix.
\end{example}
\begin{example}
  Denote by \(\mathbb{R}[x]\) the space of polynomials in a single variable \(x\), and for polynomials \(\mathbf{p}(x)\) and \(\mathbf{q}(x)\) define the multiplication rule by
  \[
    (\mathbf{p}\cdot\mathbf{q})(x)\coloneqq\mathbf{p}(x)\mathbf{q}(x).
  \]
  It is clear that this multiplication is bilinear in \(\mathbf{p}\) and \(\mathbf{q}\) and satisfies the associativity condition.
  The unit for this product is the constant polynomial \(1_A(x)=1\).
  For instance
  \[
    (x^3+1)\cdot(x^2+x) = x^5 + x^4 + x^2 + x.
  \]
  In fact, since the monomials \(\{x^n:n\ge 1\}\) form a linear basis for \(A\), the four terms on the right-hand side correspond simply to \(x^3\cdot x^2\), \(x^3\cdot x\) and so on, where we use the bilinearity of the product.
  This is an example of a commutative algebra.
\end{example}

\begin{example}
    Many structures with which one is already very familiar are just algebras in disguise, Table \ref{tab:algebras} gives a few examples.
    % , but you will find them increasingly often in the wild.

\begin{table}[h!]
\centering
\begin{tabular}{c|c|c|c|c}
% \hline
\textbf{Vector Space} & \textbf{Bilinear Operator} & \textbf{Associative} & \textbf{Commutative} & \textbf{Unitary} \\
\hline
\(\C\) & Complex product & Yes & Yes & Yes \\
\hline
$\mathbb{R}^3$ & Vector cross product: $\vec{a} \times \vec{b}$ & No & No & No \\
\hline
$\mathbb{R}[x]$ & Multiplication & Yes & Yes & Yes \\
\hline
$\mathcal{M}_{n\times n} (\R)$ & Matrix multiplication & Yes & No & Yes \\
\hline
\end{tabular}
\caption{Various examples of algebras. Some authors use the term algebra to refer to a vector space equipped with \emph{any} bilinear operation, not necessarily associative. In that sense, the cross product is an algebra that is not associative.}
\label{tab:algebras}
\end{table}
\end{example} 

The prime example of an algebra in relation to signatures is the \textit{tensor algebra}.
\begin{definition}
  Let \(V\) be a finite-dimensional vector space. The \emph{tensor algebra over \(V\)} is the vector space
  \[
    T(V)\coloneqq \bigoplus_{n\ge 0}V^{\otimes n}
    \quad \textrm{with \ \(V^{\otimes 0}\cong\mathbb{R}\mathbf{1}\).}
  \]
  The product is simply the tensor product, and its unit is the vector \(\mathbf{1}\) spanning \(V^{\otimes 0}\).
\end{definition}

The tensor algebra is therefore a \textit{graded vector space} in the sense of Definition \ref{def:gvs}, where order \(n\) tensors are placed in degree \(n\).
We stress the fact (see Definition \ref{def:inf.dsum}) that elements of \(T(V)\) are \emph{finite} sequences of tensors of arbitrary order.
For this reason, vectors in \(T(V)\) are usually called tensor (or non-commutative) polynomials.
Later on we will construct the space of tensor series, which are infinite sequences.

\begin{remark}
When \(V=\mathbb{R}^d\) the product can be written more explicitly in terms of the canonical basis \(\{\mathbf{e}_1,\dotsc,\mathbf{e}_d\}\) . 
Introducing the \textit{word notation}, recall 
\cite[Example 1.5]{chevyrev2016primer} 
\(\mathbf{e}_{i_1\dotsm i_n}\coloneqq \mathbf{e}_{i_1}\otimes\dotsm\otimes \mathbf{e}_{i_n}\in(\mathbb{R}^d)^{\otimes n}\) for \((i_1,\dotsc,i_n)\in\{1,\dotsc,d\}^n\), the product (denoted for now by $\cdot_{T(V)}$) is then defined as 
  \[
    \mathbf{e}_{i_1\dotsm i_n}\cdot_{T(V)}\mathbf{e}_{j_1\dotsm j_m} = \mathbf{e}_{i_1\dotsm i_nj_1\dotsm j_m}.
  \]
  For this reason it is commonly known as the \emph{concatenation product}.
  In this case, it corresponds to the product introduced in Definition \ref{def:tensor.prod}.
  Common notations for \(\mathbf{x}\cdot_{T(V)}\mathbf{y}\) include \(\mathbf{x}\otimes\mathbf{y}\) and \(\mathbf{x}\mathbf{y}\).
\end{remark} 

\begin{theorem}
  The tensor algebra enjoys the following universal property: given any algebra \(A\) and any \emph{linear} map \(f\colon V\to A\), there exists a unique map \(\hat{f}\colon T(V)\to A\), such that \(f(\mathbf{u}\otimes\mathbf{v})=f(\mathbf{u})\cdot_A f(\mathbf{v})\) for all \(\mathbf{u},\mathbf{v}\in T(V)\).
\end{theorem}
As is the case with the tensor product, this property actually characterizes the tensor algebra in the sense that any other algebra satisfying this property is necessarily isomorphic to \(T(V)\) for some vector space \(V\). 

We note that even though \(V\) is finite-dimensional, \(T(V)\) is always infinite dimensional since, owing to Proposition \ref{prp:tensor.basis},
\[
  \dim V^{\otimes n}=\left( \dim V \right)^n.
\]
For this reason, while the tensor algebra is a neat theoretical construction, it is not very useful for practical purposes.
There are a couple of ways of obtaining finite-dimensional versions of \(T(V)\) which preserve its structure. The most common in signature applications is \emph{truncation}.
The basic idea is that we want to preserve ``low order'' information while still retaining the algebra structure, where the meaning of ``order'' is in the sense of tensor level.
Luckily, the straightforward idea of just discarding high-order information works, with the caveat that the product has to be slightly modified.

\begin{definition}
  Given \(N\ge 1\), the \emph{level-\(N\) truncated tensor algebra} is the finite-dimensional graded vector space (recall Definition \ref{def:gvs})
  \[
    T^N(V)\coloneqq\bigoplus_{n=0}^NV^{\otimes n}
    \quad \textrm{
  % \]
  with product}\quad 
  % \[
\mathbf{x}\cdot_N\mathbf{y}=\begin{cases}\mathbf{x}\otimes\mathbf{y}&\text{ if }|\mathbf{x}|+|\mathbf{y}|\le N\\ 0&\text{ else }.\end{cases}.
  \]
\end{definition}

Following from Definition \ref{def:gvs}, we note that in particular every element of $T^N(V)$ can be written as a sequence of homogeneous elements, that is every \(\mathbf{v}\in T^N(V)\) is of the form \(\mathbf{v}=(\mathbf{v}_0,\mathbf{v}_1,\dotsc,\mathbf{v}_N) \) 
with \(\mathbf{v}_n\in V^{\otimes n}\) (with some of them eventually zero).
Hence, the product \(\cdot_N\) is well-defined for all \(\mathbf{x},\mathbf{y}\in T^N(V)\) and not just for homogeneous tensors -- elements in \(T^N(V)\) are thus finite sequence of tensors of order up to \(N\), with componentwise addition and multiplication by scalars.

It can be checked that \(T^N(V)\) is an algebra\footnote{Technically speaking \(T^N(V)\) is a quotient of \(T(V)\) by a two-sided ideal.} and
\[
  \dim T^N(V) = \frac{d^{N+1}-1}{d-1} \quad \textrm{where \ \(\dim V=d\).  }
\]

\begin{example}\label{example:product2}
  Let us take \(N=2\) and \(V=\mathbb{R}^d\). The space \(T^2(V)\cong\R\mathbf{1}\oplus\R^d\oplus(\R^d)^{\otimes 2}\) consists of elements of the form \((a,\mathbf{x},\mathbf{A})\) with \(a\in\R\), \(\mathbf{x}\in\R^d\) and \(\mathbf{A}\in\mathcal{M}_{d\times d}(\R)\) 
  % \fbox{REf03: For \(\mathbf{x}\in\R^d\) and \(\mathbf{A}\in\mathcal{M}_{d\times d}(\R)\), -- Basis fix required? Emphasise}
\footnote{
Note that we are tacitly identifying real valued 2 tensors (of shape $(2,2)$) with real valued \(2\)-by-\(2\) matrices -- looking back at Remark \ref{rem:tensor.coord} and Examples \ref{example:tensor-matrices-example} and \ref{example: tensor multiplication}, we have implicitly made the assumption of working with the canonical basis of $\R^2$.}.

  The product reads
  \[
    (a,\mathbf{x},\mathbf{A})\otimes(a',\mathbf{x}',\mathbf{A}') = (aa', a\mathbf{x}'+a'\mathbf{x}, a'\mathbf{A} + a\mathbf{A}' + \mathbf{x}\otimes \mathbf{x}').
  \]
  We remark that this product is not commutative, meaning that in general the above expression will be different from that of \((a',\mathbf{x}',\mathbf{A}')\otimes(a,\mathbf{x},\mathbf{A})\).
\end{example}

\begin{exercise}\label{exercise:tensorinverse}
  Let \(N=2\) and \(V=\mathbb{R}^d\). Show that an element \((a,\mathbf{x},\mathbf{A})\in T^2(V)\) is invertible if and only if \(a\neq 0\), and compute its inverse.
\end{exercise}
 
\begin{definition}
  The \emph{extended tensor algebra} is the direct product  
  \[
    T\dlpar V\drpar \coloneqq \prod_{n=0}^{\infty}V^{\otimes n}.
  \]
\end{definition}
We identify \(T\dlpar V\drpar\) with the space of infinite sequences \(\mathbf{u}=(\mathbf{u}_0,\mathbf{u}_1,\dotsc)\) with \(\mathbf{u}_0\in\R\), \(\mathbf{u}_1\in V\), and so on.
The product is induced by the product on \(T(V)\) and is given, for \(\mathbf{u}=(\mathbf{u}_0,\mathbf{u}_1,\dotsc)\) and \(\mathbf{v}=(\mathbf{v}_0,\mathbf{v}_1,\dotsc)\), by \(\mathbf{u}\mathbf{v}=\mathbf{w}=(\mathbf{w}_0,\mathbf{w}_1,\dotsc)\) where
\[
  \mathbf{w}_n = \sum_{k=0}^{n}\mathbf{u}_k\otimes\mathbf{v}_{n-k}\in(V^*)^{\otimes n}.
\]
This product mimics polynomial multiplication and is sometimes called the \emph{Cauchy product} for this reason.
Since this space contains arbitrarily long sequences of tensors, its elements are commonly called \emph{tensor series}.
We note that the tensor algebra \(T(V)\) is a strict subspace of \(T\dlpar V\drpar\).

For each integer \(N\ge 1\) there is a canonical projection \(\pi_N\colon T\dlpar V\drpar\to T^N(V)\), preserving multiplication, given simply by discarding tensors of degree greater than \(N\), that is,
\[
  \pi_N(\mathbf{u}_0,\mathbf{u}_1,\dotsc) = (\mathbf{u}_0,\mathbf{u}_1,\dotsc,\mathbf{u}_N).
\]
In the realm of signatures, this projection is used to produce finite-dimensional versions of the signature (see Example 3.12 just below) that are suitable for its representation in a computer.

\begin{example}
  Recall Example \ref{example:basicSigstotenosor}. Our prime example of an element in \(T\dlpar \R^d\drpar\) is the signature of a smooth \(\R^d\)-valued path \(\mathbf{x}\colon[0,1]\to\R^d\).
  Its signature over the interval \([s,t]\subseteq[0,1]\), denoted by \(S(\mathbf{x})_{s,t}\), is the tensor series of iterated integrals:
  \[
    S(\mathbf{x})_{s,t}\coloneqq\left( 1,\int_s^t\mathrm{d}\mathbf{x}_u,\int_s^t\int_s^{u_2}\mathrm{d}\mathbf{x}_{u_1}\otimes \mathrm{d}\mathbf{x}_{u_2},\dotsc \right).
  \]
  Projecting to the level-2 truncated tensor algebra we get
  \[
    \pi_2S(\mathbf{x})_{s,t}=\left( 1,\int_s^t\mathrm{d}\mathbf{x}_u, \int_s^t\int_s^{u_2}\mathrm{d}\mathbf{x}_{u_1}\otimes \mathrm{d}\mathbf{x}_{u_2}\right).
  \]
\end{example}

%%%%%%%%%%%%%%%%%%%%%%%%%%%%%%%%%%%%%%%%%%%%
%%%%%%%%%%%%%%%%%%%%%%%%%%%%%%%%%%%%%%%%%%%%
%%%%%%%%%%%%%%%%%%%%%%%%%%%%%%%%%%%%%%%%%%%%
%%%%%%%%%%%%%%%%%%%%%%%%%%%%%%%%%%%%%%%%%%%%
% \appendix
% \newpage
\section[Solutions to Exercises]{Solutions to Exercises}
\small

\begin{solution}[To Exercise \ref{ex:dsum}] 
  We must check that the operations satisfy the axioms of a vector space, that is, that \(+\) is associative and commutative, and that scalar multiplication distributes over \(+\).
  Let \((\mathbf{u}_1,\mathbf{v}_1),(\mathbf{u}_2,\mathbf{v}_2)\in U\oplus V\). Then
  \begin{align*}
    (\mathbf{u}_1,\mathbf{v}_1)+(\mathbf{u}_2,\mathbf{v}_2) &= (\mathbf{u}_1+\mathbf{u}_2, \mathbf{v}_1+\mathbf{v}_2)
    % \\
    % &
    = (\mathbf{u}_2+\mathbf{u}_1, \mathbf{v}_2+\mathbf{v}_1) 
    = (\mathbf{u}_2,\mathbf{v}_2) + (\mathbf{u}_1,\mathbf{v}_1).
  \end{align*}
  Moreover, if \((\mathbf{u}_3,\mathbf{v}_3)\in U\oplus V\) then
  \begin{align*}
    ((\mathbf{u}_1,\mathbf{v}_1)+(\mathbf{u}_2,\mathbf{v}_2))+(\mathbf{u}_3,\mathbf{v}_3) &= (\mathbf{u}_1+\mathbf{u}_2,\mathbf{v}_1+\mathbf{v}_2)+(\mathbf{u}_3,\mathbf{v}_3) 
    % \\
    % &
    = (\mathbf{u}_1+\mathbf{u}_2+\mathbf{u}_3,\mathbf{v}_1+\mathbf{v}_2+\mathbf{v}_3) 
    \\
    &= (\mathbf{u}_1,\mathbf{v}_1)+(\mathbf{u}_2+\mathbf{u}_3,\mathbf{v}_2+\mathbf{v}_3) 
    \\
    &
    = (\mathbf{u}_1,\mathbf{v}_1)+((\mathbf{u}_2,\mathbf{v}_2)+(\mathbf{u}_3,\mathbf{v}_3)).
  \end{align*} 
  Likewise, for any \(\lambda\in \mathbb{R}\) we have
  \begin{align*}
    \lambda((\mathbf{u}_1,\mathbf{v}_1)+(\mathbf{u}_2,\mathbf{v}_2))&=\lambda(\mathbf{u}_1+\mathbf{u}_2,\mathbf{v}_1+\mathbf{v}_2)\\
    &= (\lambda(\mathbf{u}_1+\mathbf{u}_2),\lambda(\mathbf{v}_1+\mathbf{v}_2)) 
    % \\
    % &
    = (\lambda \mathbf{u}_1+\lambda \mathbf{u}_2,\lambda \mathbf{v}_1+\lambda \mathbf{v}_2) \\
    &= (\lambda \mathbf{u}_1,\lambda \mathbf{v}_1) + (\lambda \mathbf{u}_2,\lambda \mathbf{v}_2) 
    = \lambda(\mathbf{u}_1,\mathbf{v}_1)+\lambda(\mathbf{u}_2,\mathbf{v}_2).
  \end{align*}
  We have used throughout that \(U\) and \(V\) are vector spaces.

  Additive inverses are given simply by \(-(\mathbf{u},\mathbf{v}) = (-\mathbf{u},-\mathbf{v})\) while the neutral element is \(0_{U\oplus V}=(0_U,0_V)\).
\end{solution}

\begin{solution}[To Exercise \ref{ex:tensor}]
  We have to check that the operations satisfy the axioms of a vector space.
  Since addition is defined symmetrically, we only check one side.
  Let \(\mathbf{u}_1,\mathbf{u}_2\in U\) and \(\mathbf{v}\in V\).
  Then
  \[
    (\mathbf{u}_1,\mathbf{v}) + (\mathbf{u}_2,\mathbf{v}) = (\mathbf{u}_1+\mathbf{u}_2,\mathbf{v}) = (\mathbf{u}_2+\mathbf{u}_1,\mathbf{v}) = (\mathbf{u}_2, \mathbf{v}) + (\mathbf{u}_1,\mathbf{v}).
  \]
  Associativity follows in a similar way:
  \begin{align*}
    \left( (\mathbf{u}_1,\mathbf{v})+(\mathbf{u}_2,\mathbf{v}) \right) + (\mathbf{u}_3,\mathbf{v}) 
&= (\mathbf{u}_1+\mathbf{u}_2,\mathbf{v})+(\mathbf{u}_3,\mathbf{v})
\\
&
= \big((\mathbf{u}_1+\mathbf{u}_2)+\mathbf{u}_3,\mathbf{v}\big)
=\big (\mathbf{u}_1+(\mathbf{u}_2+\mathbf{u}_3),\mathbf{v}\big)
\\
&= (\mathbf{u}_1,\mathbf{v})+(\mathbf{u}_2+\mathbf{u}_3,\mathbf{v}) 
% \\
% &
= (\mathbf{u}_1,\mathbf{v})+\left( (\mathbf{u}_2,\mathbf{v})+(\mathbf{u}_3,\mathbf{v}) \right).
  \end{align*}

  Now, for any \(\lambda\in\R\) we see that (using throughout that \(U\) and \(V\) are vector spaces)
  \begin{align*}
    \lambda\left( (\mathbf{u}_1,\mathbf{v})+(\mathbf{u}_2,\mathbf{v}) \right) 
    &= \lambda(\mathbf{u}_1+\mathbf{u}_2,\mathbf{v})
    % \\
% &
= (\lambda(\mathbf{u}_1+\mathbf{u}_2),\mathbf{v})
\\
&
= (\lambda\mathbf{u}_1+\lambda\mathbf{u}_2,\mathbf{v})
% \\
% &
= (\lambda\mathbf{u}_1,\mathbf{v})+(\lambda\mathbf{u}_2,\mathbf{v})
% \\
% &
= \lambda(\mathbf{u}_1,\mathbf{v}) + \lambda(\mathbf{u}_2,\mathbf{v}).
  \end{align*}

  Additive inverses are given by \(-(\mathbf{u},\mathbf{v})=(-\mathbf{u},\mathbf{v})=(\mathbf{u},-\mathbf{v})\).
\end{solution}

\begin{solution}[To Exercise \ref{ex:otimes.familiarity}]
  \begin{enumerate}[label=(\alph*)]
\item \label{item:axioms}Let \(U,V\) be vector spaces. On the set \(U\times V\) define the following bilinearity operations: for \(\mathbf{u}, \mathbf{u}_1, \mathbf{u}_2\in U\), \(\mathbf{v}, \mathbf{v}_1, \mathbf{v}_2\in V\) and \(\lambda\in\mathbb{R}\)
\begin{align*}
    (\mathbf{u}_1 \otimes \mathbf{v}) + (\mathbf{u}_2 \otimes \mathbf{v}) 
    &\coloneqq (\mathbf{u}_1 + \mathbf{u}_2) \otimes \mathbf{v}, \\
    (\mathbf{u} \otimes \mathbf{v}_1) + (\mathbf{u} \otimes \mathbf{v}_2) &\coloneqq \mathbf{u} \otimes (\mathbf{v}_1 + \mathbf{v}_2)  
    \ \textrm{ and } \ 
    % \\
    \lambda(\mathbf{u} \otimes \mathbf{v}) 
    % &
    \coloneqq (\lambda \mathbf{u}) \otimes \mathbf{v} \eqqcolon \mathbf{u} \otimes (\lambda \mathbf{v}).
\end{align*}
\item Set \(\mathbf{u}\coloneq (\mathbf{u}_1+\mathbf{u}_2)\). It is definitely not necessary to write \(\mathbf{u}\coloneq (\mathbf{u}_1+\mathbf{u}_2)\), but it may help to see how the axioms from \ref{item:axioms} can be applied. We then have, 
\begin{align*}
    (\mathbf{u}_1+\mathbf{u}_2)\otimes(\mathbf{v}_1+\mathbf{v}_2) 
    &= \mathbf{u}\otimes (\mathbf{v}_1+\mathbf{v}_2)
    % \\
    % &
    =\mathbf{u}\otimes \mathbf{v}_1 + \mathbf{u}\otimes \mathbf{v}_2
    % \\
    % &
    =(\mathbf{u}_1+\mathbf{u}_2)\otimes \mathbf{v}_1 + (\mathbf{u}_1+\mathbf{u}_2)\otimes \mathbf{v}_2
    \\
    % \Aboxed{(\mathbf{u}_1+\mathbf{u}_2)\otimes(\mathbf{v}_1+\mathbf{v}_2)
    &
    =\mathbf{u}_1\otimes \mathbf{v}_1+\mathbf{u}_2\otimes \mathbf{v}_1 + \mathbf{u}_1\otimes \mathbf{v}_1+\mathbf{u}_2\otimes \mathbf{v}_2.
    % }
\end{align*}
\item We have \(\lambda \mathbf{u}\otimes \lambda \mathbf{v} =\lambda^2 (\mathbf{u}\otimes \mathbf{v})\). 
Compare this with the linear scaling \eqref{eq:direct.sum.2} in Definition \ref{def:direct.sum}.
\end{enumerate}

\end{solution}

\begin{solution}[To Exercise \ref{ex:tensor.0}]
  It suffices to check that for any \(\mathbf{u}\in U\) and any elementary tensor \(\mathbf{u}'\otimes\mathbf{v}'\in U\otimes V\) it holds that \(\mathbf{u}\otimes 0_V + \mathbf{u}'\otimes\mathbf{v}'=\mathbf{u}'\otimes\mathbf{v}'\).

  Indeed, since we can write \(0_V=\mathbf{v}'-\mathbf{v}'\) it follows that
  \begin{align*}
    \mathbf{u}\otimes 0_V + \mathbf{u}'\otimes\mathbf{v}' 
    &= \mathbf{u}\otimes(\mathbf{v}'-\mathbf{v}')+\mathbf{u}'\otimes\mathbf{v}' 
    % \\
    % &
    = \mathbf{u}\otimes\mathbf{v}'+\mathbf{u}\otimes(-\mathbf{v}')+\mathbf{u'}\otimes\mathbf{v}' 
    % \\
    % &
    = \mathbf{u}'\otimes\mathbf{v}'
  \end{align*}
  where in the last equality we have used that \(\mathbf{u}\otimes(-\mathbf{v})\) is the additive inverse of \(\mathbf{u}\otimes\mathbf{v}\).

  The check for \(0_U\otimes\mathbf{v}\) can be done in a similar way.
\end{solution}

\begin{solution}[To Exercise \ref{ex:tensor.R}]
  We must find a bijective linear function \(\Psi\colon\R\otimes U\to U\).
  It suffices to define \(\Psi(\lambda\otimes\mathbf{u})=\lambda\mathbf{u}\).
  Linearity follows from the fact that the right-hand side is bilinear in \((\lambda,\mathbf{u})\) and the properties of the tensor product.
  Injectivity is immediate since \(\lambda\otimes\mathbf{u}\in\ker\Psi\) if and only if \(\Psi(\lambda\otimes\mathbf{u})=\lambda\mathbf{u}=0_U\), which in turn implies that either \(\lambda=0\) or \(\mathbf{u}=0_U\) by the axioms of vector spaces, and in both cases this means that \(\lambda\otimes\mathbf{u}=0_{\R\otimes U}\).
  Therefore, it follows that \(\ker\Psi=\{0\}\), i.e., \(\Psi\) is injective.

  Bijectivity can be shown by noting that every \(\mathbf{u}\in U\) can be obtained as \(\mathbf{u}=\Psi(1\otimes\mathbf{u})\) (which in particular implies that the inverse map is \(\Psi^{-1}(\mathbf{u})=1\otimes\mathbf{u}\)), or by noting that since \(\dim(\R\otimes U)=\dim(\R)\cdot\dim(U)=\dim(U)\), by the rank-nullity theorem it follows that
  \[
    \dim\big(\operatorname{im}(\Psi)\big)
    =\dim\big(\R\otimes U\big)-\dim\big(\ker(\Psi)\big)
    =\dim(U)
    \qquad
    \textrm{so that \quad \(\operatorname{im}(\Psi)=U\).}
  \]
\end{solution}

\begin{solution}[To Exercise \ref{exercise:tensorinverse}]
  The unit element in \(T^2(V)\) is \(\mathbf{1}=(1,0,0)\).
  From the product formula in Example \ref{example:product2} we see that the entries of the inverse element \((a',\mathbf{x}',\mathbf{A}')\coloneqq(a,\mathbf{x},\mathbf{A})^{-1}\) must satisfy
  \[
    aa' = 1,\quad a\mathbf{x}'+a'\mathbf{x} = 0
    \quad\textrm{and}\quad  
    \mathbf{A}+\mathbf{A}'+\mathbf{x}\otimes \mathbf{x}' = 0.
  \]
  The first equation is solvable if and only if \(a\neq 0\), in which case \(a'=a^{-1}\). Inserting this in the second equation it follows that $\mathbf{x}' = -\tfrac{1}{a^2}\mathbf{x}.$ Lastly, from the third equation we see that $\mathbf{A}' = -\tfrac{1}{a^2}\mathbf{A} + \tfrac{1}{a^3}\mathbf{x}\otimes \mathbf{x}.$ Hence, in \(T^2(V)\), we have that 
  \[
    (a,\mathbf{x},\mathbf{A})^{-1} = (a^{-1},-a^{-2}\mathbf{x}, - a^{-2}\mathbf{A} + a^{-3}\mathbf{x}\otimes \mathbf{x}).
  \]
  This can also be seen from the more general formula
  $\mathbf{A}^{-1} = \sum_{n\ge 0}(\mathbf{1}-\mathbf{A})^{\otimes n}$.
  
\end{solution}
\normalsize

%%%%%%%%%%%%%%%%%%%%%%%%%%%%%%%%%%%%%%%%%%%%
%%%%%%%%%%%%%%%%%%%%%%%%%%%%%%%%%%%%%%%%%%%%
%%%%%%%%%%%%%%%%%%%%%%%%%%%%%%%%%%%%%%%%%%%%
%%%%%%%%%%%%%%%%%%%%%%%%%%%%%%%%%%%%%%%%%%%%

%%%%%%%%%%%%%%%%%%
% Comment in or out the below lines to use the biblatex bibliography
%BEGINNING OF BIBTEX BIBLIOGRAPHY ENTRY
% \input{references}
% \bibliographystyle{plain}
\bibliographystyle{abbrv}
% \bibliography{tensorIntro.bib}

\begin{thebibliography}{1}

\bibitem{chevyrev2016primer}
I.~Chevyrev and A.~Kormilitzin.
\newblock A primer on the signature method in machine learning.
\newblock {\em arXiv preprint arXiv:1603.03788v2}, 2016.

\bibitem{Diestel2008book}
J.~Diestel, J.~H. Fourie, and J.~Swart.
\newblock {\em The metric theory of tensor products}.
\newblock American Mathematical Society, Providence, RI, 2008.
\newblock Grothendieck's r\'esum\'e{} revisited.

\bibitem{Hackbusch2019book}
W.~Hackbusch.
\newblock {\em Tensor spaces and numerical tensor calculus}, volume~56 of {\em Springer Series in Computational Mathematics}.
\newblock Springer, Cham, second edition, 2019.

\bibitem{haastad1989tensor}
J.~H{\aa}stad.
\newblock Tensor rank is {NP}-complete.
\newblock In {\em Automata, languages and programming ({S}tresa, 1989)}, volume 372 of {\em Lecture Notes in Comput. Sci.}, pages 451--460. Springer, Berlin, 1989.

\bibitem{hornMatrix}
R.~A. Horn and C.~R. Johnson.
\newblock {\em Matrix Analysis}.
\newblock Cambridge University Press, Cambridge, 2nd edition, 2012.

\bibitem{KoBa09}
T.~G. Kolda and B.~W. Bader.
\newblock Tensor decompositions and applications.
\newblock {\em SIAM Rev.}, 51(3):455--500, 2009.

\bibitem{kruskal1989rank}
J.~B. Kruskal.
\newblock Rank, decomposition, and uniqueness for {$3$}-way and {$N$}-way arrays.
\newblock In {\em Multiway data analysis ({R}ome, 1988)}, pages 7--18. North-Holland, Amsterdam, 1989.

\bibitem{rabanser2017introduction}
S.~Rabanser, O.~Shchur, and S.~G{\"u}nnemann.
\newblock Introduction to tensor decompositions and their applications in machine learning.
\newblock {\em arXiv:1711.10781}, 2017.

\bibitem{reedbarry1980methodmodernphysics}
M.~Reed and B.~Simon.
\newblock {\em Methods of modern mathematical physics. {I}}.
\newblock Academic Press, Inc. [Harcourt Brace Jovanovich, Publishers], New York, second edition, 1980.
\newblock Functional analysis.

\end{thebibliography}
%END OF BIBTEX BIBLIOGRAPHY ENTRY

%%%%%%%%%%%%%%%%%%%%%
% Comment in or out the below lines to use the biblatex bibliography
%BEGINNING OF BIBLATEX BIBLIOGRAPHY ENTRY
% \AtNextBibliography{\small} 
% \printbibliography
%END OF BIBLATEX BIBLIOGRAPHY ENTRY

%%%%%%%%%%%%%%%%%%%%%%%%%%%%%%%%%%%%%%%%%%%%
%%%%%%%%%%%%%%%%%%%%%%%%%%%%%%%%%%%%%%%%%%%%
%%%%%%%%%%%%%%%%%%%%%%%%%%%%%%%%%%%%%%%%%%%%
%%%%%%%%%%%%%%%%%%%%%%%%%%%%%%%%%%%%%%%%%%%%

% \newpage

% \input{main05-EXTRA}

\section{Extended Section: Tensor rank}
\label{sec:TensorRankSection}

The goal of Sections \ref{sec:TensorRankSection} and \ref{sec:factoringtensorProducts} is to introduce an applicable problem for the reader to grapple with: algorithmically factoring tensor product expressions to a minimal number of terms (i.e. \(u\otimes v+u\otimes w \to u\otimes(v+w)\)). While the problem itself is not critical for path signatures, through studying the problem, we will get to apply the language of tensors and tensor products to be used continuously throughout the rest of the book. 
% From another perspective, the problem of factoring tensor product expressions is perfectly self contained, and is rich and interesting in its own right, with much room for independent exploration on the part of the reader. 
For a machine learning introduction to tensors focused around the issue of tensor decomposition we refer the reader to \cite{rabanser2017introduction}. \\

Sections \ref{sec:TensorRankSection} and \ref{sec:factoringtensorProducts} also include a few programming exercises in \textit{Mathematica}. Many of these require nothing but a copy of \textit{Mathematica}, but we also make use of a few useful functions which reside in a Github repository at 
\begin{center}
\small
\url{https://github.com/jfbeda/tensor_factoring} .
\normalsize
\end{center}
From the repository, the reader can simply download and run \verb|main.nb| to access to the functions.

\subsection{Tensor rank}
\label{sec:tensor.rank}

In linear algebra, the \emph{rank} of a matrix, often defined as the number of linearly independent rows (or columns) of the matrix, is an incredibly useful property. We wish to generalize matrix rank to \emph{tensor rank}, but it is not clear how the typical definition, in terms of linearly independent rows (or columns), can be easily generalized. Instead, we begin by recharacterizing matrix rank in a way that is easily generalized to tensors of arbitrary order.

\begin{example}
\label{exa:introducing_rank}
Consider the order \(1\) tensors \(\mathbf{u}\) and \(\mathbf{v}\), and the order \(2\) tensor \(\mathbf{A}\) defined in the canonical basis (see Remark \ref{rem:tensor.coord}) by

\begin{equation}
  \mathbf{u} \coloneqq \begin{bmatrix}1\\2\end{bmatrix} \in \R^2, \quad  \mathbf{v} \coloneqq \begin{bmatrix}3\\4\end{bmatrix} \in \R^2, \quad \mathbf{A} \coloneqq \begin{bmatrix}3&4\\6&8\end{bmatrix}\in \R^2 \otimes \R^2. \label{eq:A.decomposition}
\end{equation}

Notice that \(\mathbf{A}\) may be written as a tensor product  of \(\mathbf{u}\) and \(\mathbf{v}\). That is
\begin{equation}
\mathbf{A}=\mathbf{u}\otimes \mathbf{v}, \text{ or in terms of components, }
  \mathbf{A}^{ij} = \mathbf{u}^i \mathbf{v}^j \textrm{ for }i,j=1,2.
\end{equation}
However, such a decomposition does not always exist. For example, consider \(\mathbf{B}\) given by
\begin{equation}
  \mathbf{B}\coloneqq \begin{bmatrix}1&0\\1&1\end{bmatrix} \in \R^2 \otimes \R^2.
\end{equation}

It is not possible to find order 1 tensors \(\mathbf{w},\mathbf{z}\in \R^2\) such that \(\mathbf{B}=\mathbf{w}\otimes \mathbf{z}\).
It \emph{is} however, possible to decompose \(\mathbf{B}\) into the \emph{sum} of two tensor products:

\begin{equation}
\label{eq:B.decomposition}
    \mathbf{B}
    =
    \begin{bmatrix}1&0\\1&0\end{bmatrix} 
    +
    \begin{bmatrix}0&0\\0&1\end{bmatrix} 
    = 
    \begin{bmatrix}1\\1\end{bmatrix}  
    \otimes
    \begin{bmatrix}1 \\ 0\end{bmatrix}
    + \begin{bmatrix}0 \\ 1\end{bmatrix} 
    \otimes
    \begin{bmatrix}0 \\ 1\end{bmatrix}.
\end{equation}

In fact, \textbf{the minimum number of distinct tensor products that must be summed to give a specific matrix is another way of defining matrix rank}. This number is always equivalent to the number of independent rows (or columns). Notice that \(\mathbf{A}\) has rank \(1\), and \(\mathbf{B}\) has rank \(2\) with respect to both definitions. This definition of rank is easily generalized to tensors of arbitrary order.

\end{example}

% \gbox{suggestion build a picture like
% \cite[Fig 6]{rabanser2017introduction} \url{https://arxiv.org/pdf/1711.10781} }

% \jox{It wouldn't be too hard to make something like Fig 6 in Mathcha. This visual representation seems to be pretty common, one of the other papers has it, but it doesn't do very much for me to be honest. One possible figure with actual numbers might be: https://www.mathcha.io/editor/XGxOoSr6sposMYH6YqJ6NUewnj65UdGNE65uz4ZZK0}

\begin{definition}[Tensor rank and rank decomposition]
\label{defn:rank}
Let \(V_1,\ldots,V_m\) be vector spaces.

    \begin{itemize}

      \item An order \(m\) tensor \(\mathbf{T}\in V_1 \otimes \cdots \otimes V_m\) is said to be of \emph{rank 1} if it can be written as an outer product of \(m\) order 1 tensors:  
    \begin{equation}
      \mathbf{T} = \mathbf{v}_1\otimes\dotsm\otimes\mathbf{v}_m
    \end{equation}
    for some \(\mathbf{v}_i \in V_i\).

  \item An order \(m\) tensor \(\mathbf{T}\in V_1 \otimes \cdots V_m\) is said to be of \emph{rank} \(r\) if it can be written as the sum of \(r\) rank \(1\) tensors, \textbf{and this is the smallest \(r\) for which such a decomposition is possible.} That is, \(r\) is the smallest integer such that

    \begin{equation}
    \label{aux:someq001}
      \mathbf{T} = \sum_{l = 1}^{r}\mathbf{v}_{l,1}\otimes\dotsm\otimes\mathbf{v}_{l,m}
    \end{equation}

    for some \(v_{l,i} \in V_i\). Such a decomposition of \(\mathbf{T}\) into a minimum number of rank 1 tensors is called a \emph{rank decomposition} or \emph{rank factorization}.
    \end{itemize}
\end{definition}

\begin{remark}
  For order \(2\) tensors, tensor rank is equivalent to the usual definition of matrix rank and thus provides a valid generalization. We omit a proof, but, as we saw in Example \ref{exa:introducing_rank}, \(\operatorname{Rank}(\mathbf{A})= 1\) and \(\operatorname{Rank}(\mathbf{B}) = 2\) in both the tensor rank sense of Definition \ref{defn:rank}, and in the sense of matrix rank\footnote{The astute reader may complain that we have not \emph{proven} that it is not possible to write \(\mathbf{B}\) as a single tensor product, thus in the sense of tensor rank, \eqref{eq:B.decomposition} only shows us \(\Rank (\mathbf{B}) \le 2\). Indeed, this is correct, and as we will see, while verifying an upper bound for tensor rank is easy (one must simply check the decomposition is valid), proving such a decomposition is minimal can be difficult.}. 
\end{remark}

\begin{remark}
\label{rem:basis.independent}
    Definition \ref{defn:rank} is inherently basis independent. That said, we will often fix a basis for each \(V_i\), treating the tensors as multidimensional arrays (Remark \ref{rem:tensor.coord}.)
\end{remark}

\begin{remark}[`Order', `level' and `degree' are the same, but `rank' is not.]
\label{rem:rank.v.order}
    In some contexts, particularly in physics, it is common to use the term \emph{rank} to refer to what we have called the \emph{order} of a tensor. \textbf{Do not be confused!} In our case, the \emph{order} is trivial to compute, and the \emph{rank} is difficult. The terminology itself becomes more apparent after reading Section \ref{section:AlittleBITmore} and comes from viewing the tensor algebra over \(V=\mathbb{R}^d\) as a graded vector space: \(T(V) = \bigoplus_{n\ge 0}V^{\otimes n}\). Naturally, elements of ``the homogeneous component of degree \(n\)'', \(V^{\otimes n}\), are called homogeneous tensors and \(n\) is the degree or level. 
\end{remark}

\begin{example}
    Let us verify some of the above examples in Mathematica. Mathematica has a built-in \verb|TensorProduct| function, that the reader may wish to to use. First, let us put the tensors \(\mathbf{u},\mathbf{v},\mathbf{A},\) and \(\mathbf{B}\) from Example \ref{exa:introducing_rank} into the system.

    \textbf{Input}
    \begin{lstlisting}
    u = {1,2}
    v = {3,4}
    A = {{3,4},{6,8}}
    B = {{1,0},{1,1}}
    \end{lstlisting}

    We can then indeed verify that \(\mathbf{A}=\mathbf{u}\otimes \mathbf{v}\) with:

    \textbf{Input}\footnote{\label{foot:equalequal}Recall that in many coding languages, including Mathematica, Python, and C, `\(x=y\)' sets the variable \(x\) to be \(y\), whereas `\(x==y\)' evaluates to True if \(x\) and \(y\) are equal, and False otherwise.}
    \begin{lstlisting}
    TensorProduct[u, v] == A
    \end{lstlisting}

    \textbf{Output}
    \begin{lstlisting}
    True
    \end{lstlisting}

    We have just shown that  \(\Rank(\mathbf{A})=1\) by using the tensor product definition of rank. To compute the rank of \(\mathbf{B}\), which we expect to be \(2\), we might be tempted to use Mathematica's in-built \verb|TensorRank| function. However, just as we warned in Remark \ref{rem:rank.v.order}, this function in fact computes the \emph{order} of the tensor it is fed, \textbf{not} the \emph{rank}. To compute the tensor rank of Definition \ref{defn:rank}, we could of course just use the inbuilt Mathematica function \verb|MatrixRank| because as we have seen, tensor rank and matrix rank are the same for tensors of order 2. Unsurprisingly, \verb|MatrixRank[A]| and \verb|MatrixRank[B]| return \verb|1| and \verb|2| respectively.\\

    If we wish to compute a \emph{rank decomposition} however, Mathematica does not have an in-built function. Instead, we turn to one of the functions in the provided git repository: e.g. \verb|displayRankDecompositionOrderTwo|. This function computes the rank decomposition of an order 2 tensor using a procedure that will be illustrated in Section \ref{sec:algorithms} (Algorithm \ref{algorithm:rref}). 

    \textbf{Input}
    \begin{lstlisting}
    displayRankDecompositionOrderTwo[A]
    \end{lstlisting}
    \textbf{Output}

    \begin{tcolorbox}[
    colback=white,         % Background color
    colframe=black,        % Frame color
    boxrule=0.15mm,        % Thickness of the frame
    arc=0mm,               % Square corners
    width=\linewidth,   % Make the box wider than the text width
    ]
    \begin{equation*}
        \left(
\begin{array}{cc}
 3 & 4 \\
 6 & 8 \\
\end{array}
\right)=\left(
\begin{array}{c}
 3 \\
 6 \\
\end{array}
\right)\otimes \left(
\begin{array}{c}
 1 \\
 \frac{4}{3} \\
\end{array}
\right)
\label{eq:mathematicaA}
    \end{equation*}
    \end{tcolorbox}

\textbf{Input}
    \begin{lstlisting}
    displayRankDecompositionOrderTwo[B]
    \end{lstlisting}

    \textbf{Output}
    \begin{tcolorbox}[
    colback=white,         % Background color
    colframe=black,        % Frame color
    boxrule=0.15mm,        % Thickness of the frame
    arc=0mm,               % Square corners
    width=\linewidth,   % Make the box wider than the text width
    ]
    \begin{equation*}
    \label{eq:mathematicaB}
        \left(
\begin{array}{cc}
 1 & 0 \\
 1 & 1 \\
\end{array}
\right)=\left(
\begin{array}{c}
 1 \\
 1 \\
\end{array}
\right)\otimes \left(
\begin{array}{c}
 1 \\
 0 \\
\end{array}
\right)+\left(
\begin{array}{c}
 0 \\
 1 \\
\end{array}
\right)\otimes \left(
\begin{array}{c}
 0 \\
 1 \\
\end{array}
\right)
    \end{equation*}
    \end{tcolorbox}
    
We notice that, while the rank decomposition Mathematica gives in \eqref{eq:mathematicaB} is the same as the one we gave in \eqref{eq:B.decomposition}, the decomposition it found for \(\mathbf{A}\) in \eqref{eq:mathematicaA} is not the one we gave in \eqref{eq:A.decomposition}. This illustrates the general point that rank decompositions are not unique.

\end{example}

\subsection{Computing the rank decomposition of order 2 tensors}

Section \ref{sec:tensor.rank} provides a formal definition of tensor rank, however, in terms of practically computing a rank decomposition, it is useful to work with a slightly different formalism. As we will see, finding rank decompositions for tensors of order \(2\) is straightforward, but for tensors of order \(2\) or greater it is much more difficult, in fact, it is NP-hard. In this section we explore the question of algorithmically computing rank decompositions of order 2 tensors (viewed as matrices in a particular basis), and in Section \ref{sec:surprising_properties}, we see how the problem becomes much more challenging when we consider tensors of arbitrary order.\\

Suppose \(\mathbf{M}\in \R^n \otimes \R^m\) is an order \(2\) tensor of shape \((n,m)\). By definition, we may write the components of \(\mathbf{M}\) in any basis as
    \begin{equation}
        \mathbf{M}^{ij} 
        = \sum_{l=1}^{r} \mathbf{u}_{l,1}^i \mathbf{u}_{l,2}^j. 
    \end{equation}

    This equation is simply \eqref{aux:someq001} but specified to the case of \(m=2\), and put into in coordinate form. For computational purposes, it is now useful to abandon the view of decomposing \(\mathbf{M}\) into tensor products of order one tensors \(\mathbf{u}_{l,k}\) by interpreting the components of the vectors \(\mathbf{u}_{l,1}^i\) and \(\mathbf{u}_{l,2}^j\) as components of matrices \(\mathbf{D}_{1}^{l i}\) and \(\mathbf{D}_{2}^{l j}\) respectively. This lets us write the rank decomposition of \(\mathbf{M}^{ij}\) as
    \begin{equation}
         \mathbf{M}^{ij} 
         = \sumfromto{l}{1}{r} \mathbf{D}_{1}^{l i} \mathbf{D}_{2}^{l j}.
        \label{eq:looks.like.matrix.product}
    \end{equation}
    Now, the entries of the two matrices \(\mathbf{D}_{1}^{l i}\) and \(\mathbf{D}_{2}^{lj}\) encode the rank decomposition of \(\mathbf{M}\). We now notice that \eqref{eq:looks.like.matrix.product} looks just like the matrix product in component form. That is, interpreting \(\underset{n\times m}{\mathbf{M}}\)\footnote{We use the shorthand \(\underset{n\times m}{\mathbf{M}}\) for \emph{the \(n\)-by-\(m\) matrix \(\mathbf{M}\)}.} as an \(n\)-by-\(m\) matrix of its coefficients in a particular basis, we have that \(\Rank{(\mathbf{M})} = r\) is the smallest integer such that there exist matrices \(\underset{r\times n}{{\mathbf{D}_{1}}}\) and \(\underset{r\times m}{{\mathbf{D}_{2}}}\) that satisfy:

    \begin{equation}
        \underset{n\times m}{\mathbf{M}} = \left(\underset{r\times n}{{\mathbf{D}_{1}}}\right)^T \underset{r \times m}{\mathbf{D}_{2}}.\label{eq:D1D2}
    \end{equation}

    % Relabeling to \(\mathbf{A} \coloneqq (\mathbf{D}_{1})^T\) and \(\mathbf{B} \coloneqq \mathbf{D}_{2}\), we have that \(\operatorname{Rank}(\mathbf{M}) = r\) is the smallest number such that we may write:

    % \begin{equation}
    %     \underset{n\times m}{\mathbf{M}} = \left(\underset{n\times r}{\mathbf{A}}\right)\left( \underset{r\times m}{\mathbf{B}}\right). \label{eq:aandb}
    % \end{equation}

Indeed, it is perfectly valid to define matrix rank in terms of \eqref{eq:D1D2} (e.g. see \cite[p. 13]{hornMatrix}). We see furthermore that finding \(\mathbf{D}_1\) and \(\mathbf{D}_2\) for a given \(\mathbf{M}\) amounts to finding a rank decomposition of \(\mathbf{M}\) because the coordinates of \(\mathbf{u}_{l,1}\) and \(\mathbf{u}_{l,2}\) can be recovered via

\begin{equation}
    \mathbf{u}_{l,1}^i = \mathbf{D}_1^{li}, \quad \text{and}\quad \mathbf{u}_{l,2}^j = \mathbf{D}_2^{lj}.
\end{equation}

\begin{remark}
    It is clear from \eqref{eq:D1D2} that a rank decomposition of a matrix is not unique: given a rank decomposition \(\mathbf{M}= \mathbf{D}_1^T\mathbf{D}_2\), then \(\mathbf{M} = (\mathbf{P}^T\mathbf{D}_1)^T (\mathbf{P}^{-1} \mathbf{D}_2)\) is also a valid rank decomposition for any invertible \(r\)-by-\(r\)  matrix \(\mathbf{P}\). Stronger uniqueness conditions exist for tensors of higher rank \cite{kruskal1989rank}.
\end{remark}

\subsection{Algorithms for computing the rank decomposition of a matrix}
\label{sec:algorithms}

There are fast (i.e. polynomial time) algorithms for computing \(\mathbf{D}_1\) and \(\mathbf{D}_2\) from \eqref{eq:D1D2} for an \(n\)-by-\(m\) matrix \(\mathbf{M}\). One such algorithm is detailed below using the reduced row echelon form (RREF) of \(\mathbf{M}\)\footnote{Recall the reduced row echelon form of a matrix is the result of applying elementary row operations on a matrix to bring it into a reduced form. For example, the following matrices are all in reduced row echelon form:

\begin{equation}
   \begin{bmatrix}
1 & 7 & 0 & 2 \\
0 & 0 & 1 & 2 \\
0 & 0 & 0 & 0
\end{bmatrix}, \quad \begin{bmatrix}
1 & 0 \\
0 & 1
\end{bmatrix},\quad \text{and} \quad \begin{bmatrix}
1 & 0 & 2 \\
0 & 1 & 3 \\
0 & 0 & 0
\end{bmatrix}.
\end{equation}}:

\begin{algorithm}
\caption{Rank decomposition of a matrix using the reduced row echelon form}\label{algorithm:rref}
\begin{algorithmic}[1]
\State $\mathbf{M} \gets \text{some matrix}$
\State $\mathbf{R} \gets \operatorname{RREF}(\mathbf{M})$
\Comment Reduced row echelon form
% \State $\text{pivot\_columns} = \text{sdf}$
\State $\mathbf{D}_1^T \gets \text{all columns of \(\mathbf{M}\) that are pivot columns of \(\mathbf{R}\)}$
\State $\mathbf{D}_2 \gets \text{all non-zero rows of \(\mathbf{R}\)}$
\State \textbf{assert \(\mathbf{M} == \mathbf{D}_1^T \mathbf{D_2} \)}
\State \textbf{return} \(\mathbf{D_1}\) and \(\mathbf{D_2}\)
\end{algorithmic}
\end{algorithm}

Another method, computationally slower, uses the singular value decomposition (SVD)\footnote{The SVD of a matrix \(\underset{n\times m}{\mathbf{M}}\) is a factorization to \(\mathbf{M} = \mathbf{U} \mathbf{\Sigma} \mathbf{V}^T\) where \(\underset{n \times n}{\mathbf{U}}\) and \(\underset{m\times m}{\mathbf{V}}\) are orthogonal matrices and \(\underset{n \times m}{\mathbf{\Sigma}}\) is a diagonal matrix with \(\operatorname{Rank}(\mathbf{M})\) non-zero entries.}:

\begin{algorithm}

\caption{Rank decomposition of a matrix using singular value decomposition}\label{algorithm:svd}
\begin{algorithmic}[1]
\State $\mathbf{M} \gets \text{some matrix}$
\State $\mathbf{U}, \mathbf{\Sigma}, \mathbf{V}^T \gets \operatorname{SVD}(\mathbf{M})$
\Comment Singular value decomposition
\State $r \gets \operatorname{Rank}(\mathbf{M})$ \Comment Easily acquired from the SVD
% \State $\text{pivot\_columns} = \text{sdf}$
\State $\mathbf{U}' \gets $ first \(r\) columns of \(\mathbf{U}\)
\State $\mathbf{\Sigma}' \gets $ top-left \(r\times r\) diagonal block of \(\mathbf{\Sigma}\)
\State $V' \gets $ first \(r\) columns of \(\mathbf{V}\)
\State \textbf{assert} \(C == \mathbf{U}' \mathbf{\Sigma}' \mathbf{V}'^T\)
\State \(\mathbf{D}_1^T \gets \mathbf{U}'\)
\State \(\mathbf{D_2} \gets \mathbf{\Sigma}' \mathbf{V}'^T\)
\State \textbf{assert} \(\mathbf{M} == \mathbf{D}_1^T\mathbf{D_2}\)
\State \textbf{return} \(\mathbf{D}_1\) and \(\mathbf{D}_2\)
\end{algorithmic}
\end{algorithm}

% Algorithm \ref{algorithm:svd} is much slower than Algorithm \ref{algorithm:rref}, especially when computed symbolically.

% \jox{Occasionally the algorithm placements float away from where they are supposed to be and I haven't fixed this}

Implementations of both Algorithm \ref{algorithm:rref} and \ref{algorithm:svd} in Mathematica are available through the provided git repository via \verb|rankDecompositionOrderTwoToMatricesRREF| and \verb|rankDecompositionOrderTwoToMatricesSVD| respectively.

\begin{exercise}
\label{exercise:Compute me a rank}

Define $\mathbf{M}$ as 
\[\mathbf{M} = \begin{bmatrix}
            3 & 4& 2\\
            1 & 2 & 1\\
            0 & -2 & -1\\
        \end{bmatrix}.\]

    \begin{enumerate}[label=(\alph*)]
\item Compute the rank of \(\mathbf{M}\) by considering the number of independent columns. \label{item:column}
\item Compute the rank of \(\mathbf{M}\) by considering the number of independent rows.\label{item:row}
\item Does the decomposition of \(\mathbf{M}\) into the form of \eqref{eq:valid} below contradict your answers to \ref{item:column} and \ref{item:row}?

\begin{equation}
    \begin{bmatrix}
        4 \\
        2\\
        -1
    \end{bmatrix}
    \otimes
    \begin{bmatrix}
        1 \\
        0\\
        0
    \end{bmatrix}
    + 
    \begin{bmatrix}
        1\\
        1\\
        -1
    \end{bmatrix}
    \otimes 
    \begin{bmatrix}
        -1\\
        2\\
        1
    \end{bmatrix} 
    + 
    \begin{bmatrix}
        1\\
        0\\
        0
    \end{bmatrix}
    \otimes 
    \begin{bmatrix}
        0\\
        2\\
        1
      \end{bmatrix}  
      = 
      \mathbf{M}. \label{eq:valid}
\end{equation}

\end{enumerate}

    % \end{svgraybox}

\end{exercise}

\subsection{Surprising properties of the tensor decomposition}
\label{sec:surprising_properties}

Tensor rank is filled with surprising properties, not expected from a generalization of matrix rank. 
\begin{description}
        \item[$\triangleright$] \emph{Computing tensor decompositions of tensors of order greater than \(2\) is NP-hard} \cite{haastad1989tensor}.

        This is to say that, assuming \(P \ne NP\), there is no polynomial time algorithm for computing the tensor rank of tensors of order greater than 2. For the case of order 2, the algorithms in Section \ref{sec:algorithms} are both polynomial time.
        
        \item[$\triangleright$] \emph{The rank of a real tensor can differ over \(\R\) and \(\C\)} \cite[p. 10]{kruskal1989rank}.

        We can show this through the following illustrating example from Kruskal {\cite[p. 464]{KoBa09}}.

        % \item[$\triangleright$] \emph{The probability that a randomly selected order \(3\) tensor has rank \(3\) is neither 0 nor 1 \cite{bergqvist2013exact}}.

        % This is surprising because for matrices, the probability that a randomly selected \(n\)-by-\(n\) matrix is of rank \(n\) is \(1\). Defining \emph{almost all} to mean  \emph{all except a subset of measure zero}, we mean that ``almost all \(n\)-by-\(n\) matrices have rank \(n\)". To see this for the \(2\)-by-\(2\) case, consider choosing two random vectors in the plane: the probability that both vectors are parallel is \(0\). Thus \emph{almost all} \(2 \times 2\) matrices have rank \(2\). This is not the case for higher order tensors. Amazingly, about \(71\%\) \((\approx \pi/4)\) of \(2 \times 2 \times 2\) tensors have rank \(2\), and about \(29\%\) \((\approx1 - \pi/4)\) have rank \(3\) \cite{bergqvist2013exact}.
    \end{description}

% We claimed that the rank of a real tensor can differ over \(\R\) and \(\C\). This is not true for order \(2\) tensors (matrices). We can show this through an illustrating example from Kruskal {\cite[p. 464]{KoBa09}}.
\begin{example}
\label{exa:r_vs_c}
Let \(\mathbf{Z}\in \R^2 \otimes \R^2 \otimes \R^2\) be the order \(3\) tensor of shape \((2,2,2)\) defined in the canonical basis by
    \begin{equation}
      \mathbf{Z}^{ij1} = \begin{bmatrix}1&0\\0&1\end{bmatrix}^{ij}, \quad \mathbf{Z}^{ij2} = \begin{bmatrix}0&1\\-1&0\end{bmatrix}^{ij}.
    \end{equation}

    Define \(\mathbf{u}_{l,i}\) by

\begin{alignat}{3}
    \mathbf{u}_{1,1} &= \begin{bmatrix} 1 \\ 0 \end{bmatrix}, &\quad 
    \mathbf{u}_{1,2} &= \begin{bmatrix} 1 \\ 0 \end{bmatrix}, &\quad 
    \mathbf{u}_{1,3} &= \begin{bmatrix} 1 \\ -1 \end{bmatrix}, \\
    \mathbf{u}_{2,1} &= \begin{bmatrix} 0 \\ 1 \end{bmatrix}, &\quad 
    \mathbf{u}_{2,2} &= \begin{bmatrix} 0 \\ 1 \end{bmatrix}, &\quad 
    \mathbf{u}_{2,3} &= \begin{bmatrix} 1 \\ 1 \end{bmatrix}, \\
    \mathbf{u}_{3,1} &= \begin{bmatrix} 1 \\ -1 \end{bmatrix}, &\quad 
    \mathbf{u}_{3,2} &= \begin{bmatrix} 1 \\ 1 \end{bmatrix}, &\quad 
    \mathbf{u}_{3,3} &= \begin{bmatrix} 0 \\ 1 \end{bmatrix}.
\end{alignat}

Define \(\mathbf{v}_{l,j}\) by

\begin{alignat}{3}
    \mathbf{v}_{1,1} &= \frac{1}{\sqrt{2}} \begin{bmatrix} 1 \\ -i \end{bmatrix}, &\quad 
    \mathbf{v}_{1,2} &= \frac{1}{\sqrt{2}} \begin{bmatrix} 1 \\ i \end{bmatrix}, &\quad 
    \mathbf{v}_{1,3} &= \begin{bmatrix} 1 \\ -i \end{bmatrix}, \\
    \mathbf{v}_{2,1} &= \frac{1}{\sqrt{2}} \begin{bmatrix} 1 \\ i \end{bmatrix}, &\quad 
    \mathbf{v}_{2,2} &= \frac{1}{\sqrt{2}} \begin{bmatrix} 1 \\ -i \end{bmatrix}, &\quad 
    \mathbf{v}_{2,3} &= \begin{bmatrix} 1 \\ i \end{bmatrix}.
\end{alignat}

We may decompose \(\mathbf{Z}\) as 
    \begin{equation}
    \label{eq:3terms}
        \mathbf{Z} = \mathbf{u}_{1,1}\otimes\mathbf{u}_{1,2}\otimes\mathbf{u}_{1,3}+\mathbf{u}_{2,1}\otimes\mathbf{u}_{2,2}\otimes\mathbf{u}_{2,3}+\mathbf{u}_{3,1}\otimes\mathbf{u}_{3,2}\otimes\mathbf{u}_{3,3}
    \end{equation}

Which is a decomposition in three terms, using only real coefficient. However, despite the fact that \(\mathbf{Z}\) is real, we may also decompose \(\mathbf{Z}\) as:

\begin{equation}
\label{eq:2terms}
    \mathbf{Z} = \mathbf{v}_{1,1}\otimes\mathbf{v}_{1,2}\otimes\mathbf{v}_{1,3}+\mathbf{v}_{2,1}\otimes\mathbf{v}_{2,2}\otimes\mathbf{v}_{2,3}
\end{equation}

using only two terms. The price we pay for the two term decomposition is that the coefficients of the tensors in which we expand \(\mathbf{Z}\) are complex. In this way we see clearly that the rank of a tensor depends on whether or not it is decomposed over \(\R\) or over \(\C\). The astute reader will indeed point out that we have not \emph{proven} that no two-term decomposition of \(\mathbf{Z}\) over \(\R\) exists, but it can be shown that this is the case. This example serves only as an illustration.

\begin{remark}
    The tensors \(\mathbf{Z}^{ijk}\), \(\mathbf{u}_{l,i}\) and \(\mathbf{v}_{l,j}\) above are hard coded in the git repository, and can be accessed via \verb|Z[[i,j,k]]|, \verb|u[[k,i]]| and \verb|v[[k,j]]| respectively. We can then verify \eqref{eq:3terms} and \eqref{eq:2terms} by running (recall footnote \ref{foot:equalequal} for the `==' notation):

    % \small %%%HERE
    \textbf{Input}
    \begin{lstlisting}
Z==Sum[TensorProduct[u[[k,1]],u[[k,2]],u[[k,3]]],{k,3}]
Z==Sum[TensorProduct[v[[k,1]],v[[k,2]],v[[k,3]]],{k,2}]
    \end{lstlisting}  
    \textbf{Output}
    \begin{lstlisting}
    True
    True
    \end{lstlisting}
    % \normalsize %%% HERE  
\end{remark}

\end{example}

%%%%%%%%%%%%%%%%%%%%%%%%%%%%%%%%%%%%%%%%%%%%%%%%%%%%%%%%
%%%%%%%%%%%%%%%%%%%%%%%%%%%%%%%%%%%%%%%%%%%%%%%%%%%%%%%%
%%%%%%%%%%%%%%%%%%%%%%%%%%%%%%%%%%%%%%%%%%%%%%%%%%%%%%%%
%%%%%%%%%%%%%%%%%%%%%%%%%%%%%%%%%%%%%%%%%%%%%%%%%%%%%%%%
%%%%%%%%%%%%%%%%%%%%%%%%%%%%%%%%%%%%%%%%%%%%%%%%%%%%%%%%
\section{Extended Section: Factoring tensor product expressions}
\label{sec:factoringtensorProducts}

In this section, we introduce an elegant toy problem of factoring tensor product expressions to increase our familiarity with tensors and tensor products. 
The problem is interesting in its own right and is relevant to the path signature theory but not needed in the \emph{Signatures Methods in Finance} volume. 

The problem is to \emph{find an algorithm to factor tensor product $\otimes$ expressions to a minimal number of terms}. As is often the case, this problem can be elegantly recast in the language of tensors: computing a minimal factorization amounts to finding a rank decomposition.

% There are many different definitions of the \emph{tensor product} \gbox{not sure if there are}. We have already seen one in Definition \ref{def: tensor product}, where we adopted the symbol, \(\otimes\), to refer to a product between compatibly sized tensors. Now, we will introduce another tensor product, between elements of arbitrary vector spaces. In order to distinguish from the tensor product of earlier, we will denote this by \(\ootimes\). In this section, we shall mean \(\ootimes\) when we refer to the \emph{tensor product}.\\
% \gbox{these are the same. jsut different spaces in their abstract definition.}

\subsection{The main problem}
% The problem we would like to chew on is the following: 
\emph{Given a tensor product expression comprising many terms, what is the minimum number of terms in which it may be factored using the linearity rules of the tensor product?} 

Let us consider a few simple examples. Suppose \(\mathbf{a}_1,\mathbf{a}_2,\mathbf{a}_3 \in U\), and \(\mathbf{b}_1,\mathbf{b}_2,\mathbf{b}_3 \in V\) for two vector spaces \(U\) and \(V\). Consider 
    \begin{equation} \label{eq:x1}
      \begin{split}
        \mathbf{X}_0 &\coloneqq \mathbf{a}_1 \otimes \mathbf{b}_1 + \mathbf{a}_1 \otimes \mathbf{b}_2 + \mathbf{a}_2 \otimes \mathbf{b}_ 1 + \mathbf{a}_2 \otimes \mathbf{b}_2\to(\mathbf{a}_1 + \mathbf{a}_2) \otimes (\mathbf{b}_1 + \mathbf{b}_2) % \label{eq:x0}
        \\
        \mathbf{X}_1 &\coloneqq \mathbf{a}_1 \otimes \mathbf{b}_1 +\mathbf{a}_1 \otimes \mathbf{b}_3 +  \mathbf{a}_2 \otimes \mathbf{b}_2  +\mathbf{a}_2 \otimes \mathbf{b}_3\to\mathbf{a}_1 \otimes (\mathbf{b}_1 + \mathbf{b}_3) + \mathbf{a}_2 \otimes (\mathbf{b}_2 + \mathbf{b}_3 )
      \end{split}
    \end{equation}

    We see that \(\mathbf{X}_0\), originally comprising \(4\) terms, can be factored into just \(1\) term. Whereas \(\mathbf{X}_1\), also comprising \(4\) terms can only be factored into a minimum of \(2\) terms. It is not hard to convince yourself that there is no way to factor \(\mathbf{X}_1\) into just a single term, but we would like to be able to \emph{prove} this. That is, given a tensor product expression \(\mathbf{X}\) we would like a systematic way of finding a factorization of \(\mathbf{X}\) in a minimal number of terms. As it turns out, the answer to this question lies in computing the rank decomposition of \(\mathbf{X}\). Let us start by formalizing what `minimal factorization' means to then state the problem appropriately (see Problem \ref{problem:main_problem}).
\begin{definition}[Minimal factorization]
    Let \(U,V\) be vector spaces, and let \(\mathbf{X}\in U\otimes V\) be a tensor product expression. A factorization of \(\mathbf{X}\) is called \emph{minimal} (or `\textit{minimal factorization}') if there are no factorizations of \(\mathbf{X}\) with fewer terms.
\end{definition}
We note immediately that 
    minimal factorizations are not unique. Observe that $\mathbf{X}_1$ from \eqref{eq:x1} can also be minimally factored to two terms as
\begin{equation*}
    \mathbf{X}_1 = \frac12(\mathbf{a}_1+\mathbf{a}_2)\otimes(\mathbf{b}_1+\mathbf{b}_2+2\mathbf{b}_3)+\frac12(\mathbf{a}_2-\mathbf{a}_1)\otimes(\mathbf{b}_2-\mathbf{b}_1).
\end{equation*}

\begin{problem}[Minimally factoring tensor product expressions]
    \label{problem:main_problem}
    Let \(U,V\) be vector spaces of dimension \(d_1\) and \(d_2\) with bases \(\set{\mathbf{a}_i} \) and \(\set{\mathbf{b}_j}\) respectively. Given the unfactored tensor product expression:
         \begin{equation}
        \mathbf{X} \coloneqq \sumfromto{i}{1}{d_1}\sumfromto{j}{1}{d_2} \mathbf{X}^{ij} (\mathbf{a}_i \otimes \mathbf{b}_j) = \mathbf{X}^{11} \mathbf{a}_1 \otimes \mathbf{b}_1 + \cdots +  \mathbf{X}^{d_1 d_2} \mathbf{a}_{d_1} \otimes \mathbf{b}_{d_2} \in U \otimes V. \label{eq:bruh1}
    \end{equation}

    % where the \(C_{ij}\) are known.

    Determine the minimum \(r\in \N\), such that there exist  \(\mathbf{u}_l\in U\), and \(\mathbf{v}_l\in V\) such that we may write \(\mathbf{X}\) as:

    \begin{equation}
        \mathbf{X}= \sum_{l=1}^r\mathbf{u}_l \otimes \mathbf{v}_l.
        \label{eq:bruh2}
    \end{equation}
    % \(D^{(1)}_\cdot\) and \(D^{(2)}_\cdot\) such that \(\mathbf{x}\) can be factored as:
    % \begin{equation}
    %   \mathbf{x} = \sumfromto{l}{1}{r} \bracket{\sumfromto{i}{1}{n_1}\mathbf{D}^{(1)}_{l i} {\mathbf{v}_i}} \otimes \bracket{\sumfromto{j}{1}{n_2}\mathbf{D}^{(2)}_{l j} {\mathbf{u}_j}}.    \label{eq:bruh2}
    % \end{equation}
\end{problem}

\subsection{The naive algorithm}
    Already you may be looking at Problem \ref{problem:main_problem} and seeing similarities to the tensor decompositions of the previous section. For the moment however, we ignore this and try a simple algorithm for approaching the problem. Suppose the expression we need to factor is:
    \begin{equation}
        \mathbf{X}_0 = \mathbf{a}_1 \otimes \mathbf{b}_1 + \mathbf{a}_2 \otimes \mathbf{b}_2 + \mathbf{a}_1 \otimes \mathbf{b}_3 + \mathbf{a}_2 \otimes \mathbf{b}_3.
    \end{equation}

    A naive algorithm might try grouping left terms \(\bracket{\mathbf{a}\otimes \mathbf{c} +  \mathbf{b}\otimes \mathbf{c} \to (\mathbf{a}+\mathbf{b}) \otimes \mathbf{c}}\), then right terms \(\bracket{\mathbf{a}\otimes \mathbf{b} + \mathbf{a} \otimes \mathbf{c} \to \mathbf{a} \otimes (\mathbf{b}+\mathbf{c})}\), until no further terms match. Applying such an algorithm starting starting from the left would give:
    \begin{equation}
        \mathbf{X}_0 = \mathbf{a}_1 \otimes \mathbf{b}_1 + \mathbf{a}_2 \otimes \mathbf{b}_2 + (\mathbf{a}_1 + \mathbf{a}_2) \otimes \mathbf{b}_3 \qquad \text{(3 terms)} \label{eq:shit}
    \end{equation}

    Or from the right:
    \begin{equation}
        \mathbf{X}_0 = \mathbf{a}_1 \otimes (\mathbf{b}_1 + \mathbf{b}_3) + \mathbf{a}_2 \otimes (\mathbf{b}_2 + \mathbf{b}_3) \qquad \text{(2 terms)}
    \end{equation}
    It turns out the factorization with 2 terms is minimal, but we already see the limitations of this algorithm: \emph{it can get stuck}. The factorization in \eqref{eq:shit} contains no terms that can be grouped, and yet is not minimal. This algorithm is \emph{greedy}: it takes the locally optimal decision, but there is no guarantee that the final factorization it provides is minimal.

\subsection{Using the power of tensor rank}
\label{sec:using_the_tensor_rank} 
We seek a better factorization algorithm. To proceed, we take the equation for the factored form of an arbitrary tensor product expression, \eqref{eq:bruh1}, and the equation for the unfactored form, \eqref{eq:bruh2}, and set them equal:

\begin{align}
    \mathbf{X} = \sum_{i=1}^{d_1} \sum_{j=1}^{d_2}\mathbf{X}^{ij}(\mathbf{a}_i \otimes \mathbf{b}_j) = \sum_{l=1}^r \mathbf{u}_l \otimes \mathbf{v}_l.\label{eq:equal}
\end{align}

Comparing \eqref{eq:equal} with Definition \ref{defn:rank}, we notice that finding the minimal factorization of \(\mathbf{X}\), i.e. by finding \(\mathbf{u}_l\) and \(\mathbf{v}_l\), is exactly the problem of finding a rank decomposition of \(\mathbf{X}\). We saw various algorithms in Section \ref{sec:algorithms} for computing the rank decomposition of an order 2 tensor, and we can indeed use any of these to find a decomposition of \(\mathbf{X}\).

To be more explicit, we may write \(\mathbf{u}_l\) and \(\mathbf{v}_l\) in terms of the bases \(\set{\mathbf{a}_i}\) and \(\set{\mathbf{b}_i}\) respectively. That is,

\begin{equation}
        \mathbf{u}_l=\sum_{i=1}^{d_1}\mathbf{u}_l^i \mathbf{a}_i,\quad \text{and}\quad
        \mathbf{v}_l=\sum_{j=1}^{d_2}\mathbf{v}_l^j \mathbf{b}_i.
\end{equation}

Substituting this into \eqref{eq:equal} gives:

\begin{align}
    \sum_{i=1}^{d_1} \sum_{j=1}^{d_2}\mathbf{X}^{ij}(\mathbf{a}_i \otimes \mathbf{b}_j) &= \sum_{l=1}^r \sum_{i=1}^{d_1}\sum_{j=1}^{d_2} \mathbf{u}^i_l \mathbf{v}_l^j \mathbf{a}_i \otimes \mathbf{b}_j\\
    0&=\sum_{i=1}^{d_1}\sum_{j=1}^{d_2} \bracket{\mathbf{X}^{ij}-\sum_{l=1}^r \mathbf{u}_l^i \mathbf{v}_l^j} \mathbf{a}_i \otimes \mathbf{b}_j.\label{eq:coefficients}\\
    \implies\quad  \mathbf{X}^{ij} &= \sum_{l=1}^r \mathbf{u}_l^i \mathbf{v}_l^j.
\end{align}

Where in the last line we have used the fact that \(\set{\mathbf{a}_i\otimes \mathbf{b}_j}\) are a basis for \(U \otimes V\) (Proposition \ref{prp:tensor.basis}) and so the coefficients in \eqref{eq:coefficients} must all vanish. We may now use any of the algorithms in Section \ref{sec:algorithms} to find matrices \(\underset{r\times n_1}{\mathbf{D}_1}\) and \(\underset{r\times n_2}{\mathbf{D}_2}\) such that \(\mathbf{X}^{ij}= (\mathbf{D}_1^T\mathbf{D}_2)^{ij}\), which allows us to recover the coefficients of \(\mathbf{u}_l\) and \(\mathbf{v}_l\) in the \(\set{\mathbf{a}_i}\) and \(\set{\mathbf{b}_i}\) bases via

\begin{equation}
    \mathbf{u}_l^i=\mathbf{D}_1^{li}, \quad \text{and} \quad \mathbf{v}_l^j = \mathbf{D}_2^{lj}.
\end{equation}

\begin{exercise}
\label{exercise:E}
    Let \(U=\R[x]\), and \(V=\R[y]\) be the vector spaces of polynomials in the single variable \(x\), and \(y\) respectively. Consider the following expression:

    \begin{equation}
        \mathbf{E}\coloneq -x\otimes y +2x^2\otimes y +3x\otimes y^2 -4 x^2\otimes y^2 + x^3 \otimes y^2 \in U\otimes V.
        \label{eq:exercise.E}
    \end{equation}

    \begin{enumerate}[label=(\alph*)]
    \item Try to factor \eqref{eq:exercise.E} to a minimal number of terms by grouping.
    \item Confirm your factorization is indeed minimal by using Mathematica, or otherwise.
    \end{enumerate}
\end{exercise}

\subsection{Generalization}

We have solved the problem of factorizing tensor product expressions containing one tensor product symbol per term. That is, we have a computationally efficient\footnote{Since all of the algorithms in Section \ref{sec:algorithms} are polynomial time, factorizations of tensor product expressions involving a single tensor product symbol are fast to compute.} algorithm for solving Problem \ref{problem:main_problem}. We can easily generalize this to the case of expressions that contain more than one tensor product symbol in each term. For example
\begin{equation}
        \mathbf{Z} \coloneqq \mathbf{u}_1 \otimes \mathbf{v}_1 \otimes \mathbf{w}_1 + \mathbf{u}_1 \otimes \mathbf{v}_2 \otimes \mathbf{w}_2 -\mathbf{u}_2 \otimes \mathbf{v}_1 \otimes \mathbf{w}_2 + \mathbf{u}_2 \otimes \mathbf{v}_2 \otimes \mathbf{w}_1. \label{eq:y0}
 \end{equation}

% Recall that the tensor product is linear in all of its arguments:
% \begin{align}
%     \nonumber
%     \mathbf{a} \otimes\cdots \otimes
%     &(\beta \mathbf{u}_1 + \gamma \mathbf{u}_2) \otimes\cdots \otimes \mathbf{z} = \\
%     & \qquad \qquad \beta (\mathbf{a} \otimes\cdots \otimes \mathbf{u}_1 \otimes\cdots \otimes \mathbf{z}) + \gamma( \mathbf{a} \otimes \cdots \otimes \mathbf{u}_2 \otimes\cdots \otimes \mathbf{z}).
% \end{align}

Proceeding just as in Problem \ref{problem:main_problem} we can state the question of factoring tensor product expressions involving more than one tensor product symbol.

\begin{problem}
\label{problem:more_than_one}
Let \(V_1, \cdots, V_m\) be vector spaces of dimension \(d_1, \cdots, d_m\) with bases \(\set{\mathbf{a}_{1,i_1}},\cdots,\set{\mathbf{a}_{m,i_m}}\) respectively. Given the unfactored tensor product expression
\begin{equation}
        \mathbf{Y} \coloneqq \sumfromto{i_1}{1}{d_1}\cdots \sumfromto{i_m}{1}{d_m} \mathbf{Y}^{i_1 i_2 \cdots i_m} \bracket{\mathbf{a}_{1,i_1} \otimes \cdots \otimes \mathbf{a}_{m,i_m}} \in V_1 \otimes \cdots \otimes V_m \label{eq:yappy}
    \end{equation}

    what is the minimum integer \(r\) such that we may write \(\mathbf{Y}\) as

    \begin{equation}
        \mathbf{Y}= \sum_{l=1}^r \mathbf{v}_{l,1} \otimes \cdots \otimes \mathbf{v}_{l,m}.
    \end{equation}
\end{problem}

Once again, we see that finding the minimum factorization of a tensor product expression \(\mathbf{Y}\) amounts simply to finding the rank decomposition of \(\mathbf{Y}\). The number of terms of terms in a minimal factorization is simply the tensor rank.\\

It is here we can make use of some of the important properties of the tensor rank that we saw in Section \ref{sec:surprising_properties}.

\begin{itemize}
    \item Since minimally factoring a tensor product expression is equivalent to finding its rank decomposition, both problems are NP-hard for tensors of order greater than \(2\).
    \item The fact that the tensor rank of real tensors can differ over \(\R\) and over \(\C\) translates into the fact that tensor product expressions of order greater than \(2\) differ when the factorization is over \(\R\) or over \(\C\). As an example, take \(\mathbf{Z}\) from equation \eqref{eq:y0}:
    \begin{equation}
        \mathbf{Z} = \mathbf{u}_1 \otimes \mathbf{v}_1 \otimes \mathbf{w}_1 + \mathbf{u}_1 \otimes \mathbf{v}_2 \otimes \mathbf{w}_2 -\mathbf{u}_2 \otimes \mathbf{v}_1 \otimes \mathbf{w}_2 + \mathbf{u}_2 \otimes \mathbf{v}_2 \otimes \mathbf{w}_1.
    \end{equation}

    This purely real expression can be minimally \textit{factored to three terms} over \(\R\):
 \begin{equation}
        \mathbf{Z} = \mathbf{u}_1 \otimes \mathbf{v}_1 \otimes (\mathbf{w}_1 -\mathbf{w}_2) + \mathbf{u}_2 \otimes \mathbf{v}_2 \otimes (\mathbf{w}_1 + \mathbf{w}_2) + (\mathbf{u}_1 -\mathbf{u}_2) \otimes (\mathbf{v}_1 + \mathbf{v}_2) \otimes \mathbf{w}_2,
    \end{equation}
    but minimally \textit{factored to two terms} over \(\C\):
    \begin{align}
    \nonumber
        \mathbf{Z} & = \frac{1}{2}\Big((\mathbf{u}_1-i\mathbf{u}_2) \otimes (\mathbf{v}_1+i\mathbf{v}_2) \otimes (\mathbf{w}_1 -i\mathbf{w}_2)
        \\
        & \qquad \qquad \qquad +  (\mathbf{u}_1+i\mathbf{u}_2) \otimes (\mathbf{v}_1-i\mathbf{v}_2) \otimes (\mathbf{w}_1 +i\mathbf{w}_2)\Big).
    \end{align}

%     In fact, the tensor \(\mathbf{Z}\), expressed in the canonical basis, is exactly the \(2\)-by-\(2\)-by-\(2\) tensor of Example \ref{exa:r_vs_c} from Kruskal \cite{KoBa09} (Fig. \ref{fig:xijk}).

% \begin{figure}
%     \centering
%     \includegraphics[width=0.4\linewidth]{images/tensor.png}
%     \caption{The real tensor \(\mathbf{Z} \in (\R^2)^{\otimes 3}\) has rank \(3\) over \(\R\) but rank \(2\) over \(\C\).}
%     \label{fig:xijk}
% \end{figure}

\end{itemize}

%%%%%%%%%%%%%%%%%%%%%%%%%%%%%%%%%%%%%%%%%%%%%%%%%%%%%%%%
%%%%%%%%%%%%%%%%%%%%%%%%%%%%%%%%%%%%%%%%%%%%%%%%%%%%%%%%
%%%%%%%%%%%%%%%%%%%%%%%%%%%%%%%%%%%%%%%%%%%%%%%%%%%%%%%%
%%%%%%%%%%%%%%%%%%%%%%%%%%%%%%%%%%%%%%%%%%%%%%%%%%%%%%%%
%%%%%%%%%%%%%%%%%%%%%%%%%%%%%%%%%%%%%%%%%%%%%%%%%%%%%%%%

\section{Extended Section: Solutions to additional Exercises}
\small

\begin{solution}[To Exercise \ref{exercise:Compute me a rank}]
% \small
\begin{enumerate}[label=(\alph*)]
\item %The rank of \(\mathbb{M}\) is \(2\). 
Label the columns \(\mathbf{C}_1, \mathbf{C}_2, \mathbf{C}_3\). Notice that \(\mathbf{C}_2 = 2 \mathbf{C}_3\) and hence \(\mathbf{C}_2\) is linearly dependent on \(\mathbf{C}_3\). Since \(\mathbf{C}_1\) and \(\mathbf{C}_3\) are linearly independent, the matrix \(\mathbf{M}\) has a total of \(2\) linearly independent columns. Thus the rank of \(\mathbf{M}\) is \(2\).

\item Label the rows \(\mathbf{R}_1, \mathbf{R}_2, \mathbf{R}_3\). Note that \(\mathbf{R}_1 = 3 \mathbf{R}_2 + \mathbf{R}_3\), thus \(\mathbf{R}_1\) is linearly dependent on \(\mathbf{R}_2\) and \(\mathbf{R}_3\). Since \(\mathbf{R}_2\) and \(\mathbf{R}_3\) are linearly independent, the rank of \(\mathbf{M}\) is \(2\).

\item One might think that \eqref{eq:valid} provides a rank decomposition of \(\mathbf{M}\) in three terms, and thus conclude that \(\Rank(\mathbf{M})=3\), contradicting our results from \ref{item:column} and \ref{item:row}. However, just because the decomposition in \eqref{eq:valid} does correctly produce \(\mathbf{M}\), it is not a decomposition with a minimal number of terms. In fact, from the previous parts of the question we \emph{know} that \(\mathbf{M}\) can be decomposed into the sum of two outer products. Finding these can be difficult (and the reader is not expected to be able to do so in their head), but here is one such decomposition:
\begin{equation}
    \begin{bmatrix}
        1\\
        0\\
        1
    \end{bmatrix}\otimes\begin{bmatrix}
        1\\
        0\\
        0
    \end{bmatrix} + \begin{bmatrix}
        2\\
        1\\
        -1
    \end{bmatrix}\otimes \begin{bmatrix}
        1\\
        2\\
        1
      \end{bmatrix} = \mathbf{M}.
\end{equation}
We can guarantee that there is no decomposition of \(\mathbf{M}\) in fewer than \(2\) terms because we saw earlier that \(\operatorname{Rank}(\mathbf{M}) = 2\). We can also use the provided Mathematica code to find a rank decomposition of \(\mathbf{M}\).

\textbf{Input}
\begin{lstlisting}
    M = {{3,4,2},{1,2,1},{0,-2,-1}};
    displayRankDecompositionOrderTwo[M]
\end{lstlisting}
\normalsize

\textbf{Output}
\begin{tcolorbox}[
    colback=white,         % Background color
    colframe=black,        % Frame color
    boxrule=0.15mm,        % Thickness of the frame
    arc=0mm,               % Square corners
    width=\linewidth,   % Make the box wider than the text width
    ]
\begin{equation*}
\left(
\begin{array}{ccc}
 3 & 4 & 2 \\
 1 & 2 & 1 \\
 0 & -2 & -1 \\
\end{array}
\right)=\left(
\begin{array}{c}
 3 \\
 1 \\
 0 \\
\end{array}
\right)\otimes \left(
\begin{array}{c}
 1 \\
 0 \\
 0 \\
\end{array}
\right)+\left(
\begin{array}{c}
 4 \\
 2 \\
 -2 \\
\end{array}
\right)\otimes \left(
\begin{array}{c}
 0 \\
 1 \\
 \frac{1}{2} \\
\end{array}
\right)
\end{equation*}
\end{tcolorbox}

This is different, valid rank decomposition of \(\mathbf{M}\).

\end{enumerate}
% \normalsize
\end{solution}

\begin{solution}[To Exercise \ref{exercise:E}]

\begin{enumerate}[label=(\alph*)]
    \item We begin with:

    \begin{equation}
        \mathbf{E}\coloneq -x\otimes y +2x^2\otimes y +3x\otimes y^2 -4 x^2\otimes y^2 + x^3 \otimes y^2 \in U\otimes V.
    \end{equation}

    Grouping from the left gives:

    \begin{equation}
        \mathbf{E}=x\otimes(3y^2-y) + x^2 \otimes (-4y^2+2y) +x^3\otimes y^2.
    \end{equation}

    Which cannot be further grouped. Trying instead from the right gives:

    \begin{equation}
        \mathbf{E}= (-x+2x^2)\otimes y + (3x-4x^2+x^3)\otimes y^2.
    \end{equation}

    Thus the rank of \(\mathbf{E}\) is no larger than \(2\). 
    
    \item Let us assign the obvious bases \(\set{x,x^2,x^3,\cdots}\) and \(\set{y,y^2,y^3,\cdots}\) to \(U=\R[x]\) and \(V=\R[y]\). Written in this basis we have

    \begin{equation}
        \mathbf{E}=\sum_{i=1}^3 \sum_{j=1}^2 \mathbf{E}^{ij} x^i\otimes y^j.
    \end{equation}

    Where the coefficients \(\mathbf{E}^{ij}\) can be written in a matrix as:

    \begin{equation}
        \mathbf{E}^{ij} = \begin{bmatrix}
            -1 & 3\\
            2 & -4\\
            0 & 1
        \end{bmatrix}.
    \end{equation}

    Letting Mathematica compute the rank of this matrix yields \(2\), which tells us that our factorization with \(2\) terms is minimal.
    \end{enumerate}
\end{solution}

% For multibibliography that we are not using
% \bibliographystyleExtra{unsrt}
% \bibliographyExtra{extraBibliography.bib}

\end{document}